\documentclass[12pt,twoside]{amsart}\usepackage{amssymb}
\usepackage{enumerate}

\newtheorem{thmspec}{\relax}
\newtheorem{theorem}{Theorem}[section]
\newtheorem{thm}[theorem]{Theorem}
\newtheorem{lem}[theorem]{Lemma}
\newtheorem{cor}[theorem]{Corollary}
\newtheorem{prop}[theorem]{Proposition}
\newtheorem{defi}[theorem]{Definition}

\theoremstyle{definition}

\theoremstyle{remark}

\numberwithin{equation}{section}
\def \Bbb{\mathbb}

\def\onto{{\kern3pt\to\kern-8pt\to\kern3pt}}

\def\<{\langle}
\def\>{\rangle}
\def\|{{\ |\ }}

\def\onto{\twoheadrightarrow}

\def\-{\underline}

\def\Re{\operatorname{Re}}
\def\Im{\operatorname{Im}}
\def\mes{\operatorname{mes}}

\def\arg{\operatorname{arg}}
\def\rad{\operatorname{rad}}
\def\Log{\operatorname{Log}}

\def\N{\Bbb N}

\def\R{\Bbb R}
\def\C{\Bbb C}

\def\B{\Bbb B}

\def\H{\Bbb H}

\def\X{\Bbb X}

\def\<{\langle}
\def\>{\rangle}

\catcode`\@=11
\def\serieslogo@{\relax}
\def\@setcopyright{\relax}
\catcode`\@=12

\title[Boundary cross theorem in dimension $1$]
{Boundary cross theorem in dimension $1$
}

\begin{document}

\author{Peter Pflug}
\address{Peter Pflug\\
Carl von  Ossietzky Universit\"{a}t Oldenburg \\
Fachbereich  Mathematik\\
Postfach 2503, D--26111\\
 Oldenburg, Germany}
\email{pflug@mathematik.uni-oldenburg.de}

\author{Vi{\^e}t-Anh  Nguy\^en}
\address{Vi{\^e}t-Anh  Nguy\^en\\
Mathematics Section\\
The Abdus Salam international centre
 for theoretical physics\\
Strada costiera, 11\\
34014 Trieste, Italy}
\email{vnguyen0@ictp.trieste.it}



\begin{abstract}
Let $X,\ Y$  be two complex manifolds of dimension $1$ which are  countable at infinity,
  let $D\subset X,$ $ G\subset Y$ be two open sets, let
  $A$ (resp. $B$) be a subset of  $\partial D$ (resp.
  $\partial G$), and let $W$ be the  $2$-fold cross $((D\cup A)\times B)\cup (A\times(B\cup G)).$
  Suppose in addition  that    $D$ (resp. $G$) is  {\it Jordan-curve-like
   on $A$} (resp. $B$) and that
  $A$ and $B$ are {\it  of positive  length}.
   We  determine the ``envelope of holomorphy"
  $\widehat{W}$ of $W$ in the sense that any function  locally bounded on $W,$
  measurable on $A\times B,$ and separately
  holomorphic
  on $(A\times G) \cup (D\times B)$ ``extends" to a function
   holomorphic  on the interior of $\widehat{W}.$
\end{abstract}
\maketitle

\section{ Introduction}

In this paper  we consider  a boundary version of the cross
theorem
in the spirit of the pioneer  work of  Malgrange--Zerner \cite{ze}.
Epstein's survey article \cite{ep} gives a historical discussion and  motivation for this
kind of theorems.

The  first results in this direction are obtained by
Komatsu \cite{ko} and Dru\.{z}kowski \cite{dr}, but only for some special
cases.
Recently, Gonchar \cite{go1,go2} has proved  a  more
general result for the one-dimensional case.
 In  recent works  \cite{pn,pn2}, the authors are able to  generalize Gonchar's result  to the higher
dimensional case.

However, in all  cases
considered so far in the literature
the hypotheses on the function to extend and its   domain of definition
are, in some sense, rather restrictive. 
 Therefore, the main goal of this work is to
establish  some boundary cross theorems in  more general
(one-dimensional) cases with more optimal hypotheses. Perhaps, this will be a first step
towards understanding the higher dimensional case in its full generality.

Our approach here  is based on the  previous work \cite{pn}, the
{\it Gonchar--Carleman operator} developed in \cite{go1,go2},  a
new result of Zeriahi \cite{zr} and a thorough geometric study of
harmonic measures.

\smallskip

\indent{\it{\bf Acknowledgment.}}   The first version of this article  was written in
Spring 2004 while the second author was visiting the Carl von
Ossietzky Universit\"{a}t Oldenburg being  supported by The
Alexander von Humboldt Foundation. He wishes to express his
gratitude to these organizations.

\section{Preliminaries}

In order to recall  the classical     versions of the boundary cross theorem and to
discuss in more detail  our motivation, we need to introduce some notation and terminology.
 In fact, we keep the main notation from the previous work \cite{pn}.
  Here $E$ denotes the open  unit disc in $\C$ and   $\mes$  the  linear measure (i.e. the one-dimensional
   Hausdorff measure).
Throughout the paper,
 for a topological space  $M,$ $\mathcal{C}(M)$ denotes the space of  all continuous functions  $f:\ M\longrightarrow\C$
equipped with the sup-norm   $\vert f\vert_M:=\sup_M \vert f\vert.$
Moreover,  a function $f:\ M\longrightarrow\C$ is said to be {\it locally bounded} on $M$ if,
for any point $z\in M,$ there are an open neighborhood $U$ of $z$ and a
 positive number $K=K_z$ such that $\vert f\vert_U <K.$ Finally, for a complex manifold $\Omega,$ $\mathcal{SH}(\Omega)$
  (resp.
$\mathcal{O}(\Omega)$) denotes the  set of
all subharmonic (resp. holomorphic) functions on $\Omega.$

In this work all complex manifolds are supposed to be countable at
infinity.

\subsection{Open set   with partly Jordan-curve-like boundary}


 Let  $X$ be a complex manifold of dimension $1.$
A {\it Jordan curve} in $X$ is the image $\mathcal{C}:=\{\gamma(t):\ t\in[a,b]\}$ of a continuous one-to-one map    $
\gamma:\ [a,b]\longrightarrow  X,$ where $a,b\in \R,\ a<b.$
The set  $ \{\gamma(t):\ t\in(a,b)\}$ is said to be the {\it interior}
of the Jordan curve.     A {\it Jordan domain}  is the image $\{\Gamma(t), \ t\in E\}$
of  a   one-to-one continuous  map    $
\Gamma:\    \overline{E}\longrightarrow  X.$
 A {\it closed Jordan curve} is the boundary  of a Jordan domain.

Consider an open set   $D\subset X.$ Then $D$   is said to be {\it
Jordan-curve-like  at a point} $\zeta\in\partial D$ if there is a
Jordan domain  $U\subset X$  such that   $\zeta\in U$ and $U\cap
\partial D$ is the interior of a  Jordan curve. Then  $\zeta$ is
said to be {\it of type 1} if there is a neighborhood $V$ of
$\zeta$ such that $V\cap D$ is a Jordan domain.  Otherwise,
$\zeta$ is said to be {\it   of type 2}. We see easily that if
$\zeta$ is of type 2, then there are an open neighborhood $V$ of
$\zeta$ and two Jordan domains $V_1,$ $V_2$  such that  $V\cap
D=V_1\cup V_2.$ Moreover, $D$  is said to be  {\it
Jordan-curve-like on a subset} $A$ of $\partial D$ if  $D$ is
Jordan-curve-like  at all points of $A.$

Now let $D\subset X$ be an open set  which is Jordan-curve-like  on a set
$A\subset\partial D.$ In the remaining part of this subsection we will introduce various
notions. We like to point out that  these notions are intrinsic, i.e., they do not depend
on any choice (of open neighborhoods, Jordan domains, conformal mappings ...)  we made in their definitions.

$A$   is said to be  {\it Jordan-measurable}
 if  for every $\zeta\in A$ the following
condition  is fulfilled:\\
{\bf Case 1: $\zeta$ is  of type 1.}
There are an open neighborhood $U=U_{\zeta}$ of $\zeta$ such that
$U\cap D$ is a   Jordan domain  and a conformal
mapping $\Phi=\Phi_{\zeta}$ from $U\cap D$ onto the unit disc $E$   which extends homeomorphically from
$\overline{U\cap D}$ onto $\overline{E}$  such that $\Phi
(\overline{U\cap D}\cap  A)$ is  Lebesgue measurable on $\partial E$.\\
{\bf Case 2: $\zeta$ is of type 2.}
There are an open neighborhood $U=U_{\zeta}$ of $\zeta$ such that $U\cap D=U_1\cup
U_2$ with Jordan  domains $U_1=U_{1,\zeta}, \  U_2= U_{2,\zeta},$ and  conformal mappings
$\Phi_j=\Phi_{j,\zeta}$ $(j=1,2)$  from  $U_j$  onto $E$ which extend homeomorphically from
$\overline{U_j}$ onto $\overline{E}$ such that
 $\Phi_j
(\overline{U_j}\cap A)$ is  Lebesgue measurable (on $\partial E$).

A  Jordan-measurable set $A\subset \partial D$ is said to be  {\it of zero length}
  if for all  $\zeta\in A,$ if one takes  $U_{\zeta},\ \Phi_{\zeta}$  when  $\zeta$ is  of type 1
  (resp. $U_{\zeta},\ \Phi_{j,\zeta}$ when $\zeta$ is of type 2)  as in the previous definition and notation, then
 $\mes\Big(\Phi_{\zeta}
(\overline{U_{\zeta}\cap D}\cap A)\Big)=0$
(resp.  $\mes\Big(\Phi_{j,\zeta}
(\overline{U_{j,\zeta}}\cap A)\Big)=0,$  $j=1,2$).

A  Jordan-measurable set $A\subset \partial D$ is said to be  {\it of positive length} if
it is not of zero length.

 Suppose that $D$ is Jordan-curve-like at a point $\zeta\in\partial D.$ We  define the
 concept of {\it angular approach regions} at  $\zeta$ as follows. 
For any $0<\alpha <\frac{\pi}{2},$   the {\it Stolz region} or
{\it angular approach region} $\mathcal{A}_{\alpha}(\zeta)$ is given by:\\
{\bf Case 1: $\zeta$ is  of type 1.}
\begin{equation*}
\mathcal{A}_{\alpha}(\zeta):=
\Phi^{-1} \left\lbrace          t\in E:\ \left\vert
 \arg\left(\frac{\Phi(\zeta)-t}{\Phi(\zeta)}\right)
 \right\vert<\alpha\right\rbrace   ,
\end{equation*}
where  $\arg:\ \C\longrightarrow (-\pi,\pi]$ is  as usual the argument function.\\
{\bf Case 2: $\zeta$ is of type 2.}
\begin{equation*}
\mathcal{A}_{\alpha}(\zeta):=\bigcup_{j=1,2}
\Phi_j^{-1} \left\lbrace          t\in E:\ \left\vert
 \arg\left(\frac{\Phi_j(\zeta)-t}{\Phi_j(\zeta)}\right)
 \right\vert<\alpha\right\rbrace  .
\end{equation*}
 Geometrically, $\mathcal{A}_{\alpha}(\zeta)$ is the
intersection of $D$ with one or  two ``cones" of aperture $2\alpha$ and vertex $\zeta$
according to the type of $\zeta.$

Let  $\zeta\in \partial D$  a point at which $D$ is Jordan-curve-like and let $U$ be an open neighborhood of
$\zeta.$
We say that  a function $f$ defined on  $U\cap D$
  admits the {\it  angular limit} $\lambda$ at $\zeta$
if
\begin{equation*}
\lim\limits_{ z\in \mathcal{A}_{\alpha}(\zeta),\ z\to\zeta} f(z)=\lambda,
\end{equation*}
for all $0<\alpha<\frac{\pi}{2}.$

Let $A\subset \partial D$ be a Jordan-measurable set and $f:\ D\longrightarrow\C,$ $g: A \longrightarrow\C$
two functions.  Then $f$ is said  to {\it  have the angular limit $g(a)$ for Jordan a.e. $a\in A,$}
if the set
\begin{equation*}
\left \lbrace a\in A:\ f\ \text{does not admit the angular limit}\ g(a)\ \text{at} \ a  \right\rbrace
\end{equation*}
is of zero length. For simplicity, in the future  we only write  ``a.e." instead of ``Jordan a.e.".

We conclude this subsection with a simple example which  may clarify  the above definitions.
Let $G$ be the open square in $\C$
whose four vertices are $1+i,$ $-1+i,$ $-1-i,$ and $1-i.$
Define the domain
\begin{equation*}
D:=G\setminus \left [-\frac{1}{2},\frac{1}{2}\right].
\end{equation*}
Then   $D$ is Jordan-curve-like on $\partial G\cup  \left
(-\frac{1}{2},\frac{1}{2}\right).$ Every point of $\partial G$ is of type
1 and every point of $\left(-\frac{1}{2},\frac{1}{2}\right)$ is of type
2.
\subsection{Harmonic measure}

Let $X$  be a complex manifold of dimension $1,$ let  $D$ be an  open subset of $X$  and let
$A\subset \partial D.$
Consider the characteristic  function
\begin{equation*}
 1_{\partial D\setminus A}(\zeta):=
\begin{cases}
1,
  &\zeta\in  \partial D\setminus A,\\
 0, & \zeta \in A.
\end{cases}
\end{equation*}
 Then  the {\it  harmonic measure} of the set
 $\partial D\setminus A$ (denoted by $\omega(\cdot,A,D)$) is the
 Perron solution of the generalized Dirichlet problem with boundary data $1_{\partial D\setminus
 A}.$ In other words, one has
\begin{equation*}
 \omega(\cdot,A,D):=\sup\limits_{u\in \mathcal{U}}u,
\end{equation*}
where $\mathcal{U}=\mathcal{U}(A,D)$ denotes the family of all subharmonic functions $u$ on
$D$ such that $\limsup\limits_{D\ni z\to\zeta} u(z)\leq  1_{\partial D\setminus A}(\zeta)$
for each $\zeta\in\partial D.$

It is well-known  (see, for example,  the book of Ransford \cite{ra} for the case $X:=\C$) that
$ \omega(\cdot,A,D)$ is harmonic on $D.$

For a point
 $\zeta\in  \partial D $ at which   $D$ is a Jordan-curve-like, we say that it is   {\it a  locally regular point relative to}
  $A$
if
\begin{equation*}
\lim\limits_{z\to \zeta,\ z\in
\mathcal{A}_{\alpha}(\zeta)}  \omega(z,A\cap U,D\cap U)=0
\end{equation*}
for any $0<\alpha <\frac{\pi}{2}$ and any open neighborhood $U$ of $\zeta.$ Obviously,
$\zeta\in \overline{A}.$
If, moreover, $\zeta\in A,$ then $\zeta$ is said to be
a {\it  locally regular point of}   $A.$
The set of all locally regular points relative to $A$ is denoted by $A^{\ast}.$
Observe that, in general, $A^{\ast}\not\subset A,$  $A\not\subset
A^{\ast}.$ However, if $A$ is open in $\partial D$ and  $D$ is  Jordan-curve-like on $A,$
 then  $A \subset A^{\ast}.$

As an immediate consequence of the Subordination Principle   for  the harmonic measure
(see Corollary 4.3.9 in  \cite{ra}), one gets
  \begin{equation}\label{eq2.1}
\lim\limits_{z\to \zeta,\ z\in
\mathcal{A}_{\alpha}(\zeta)}  \omega(z,A,D)=0,\qquad \zeta\in A^{\ast},\ 0<\alpha<\frac{\pi}{2}.
\end{equation}
We extend the function $\omega(\cdot,A,D)$ to $D\cup A^{\ast}$  by simply setting
\begin{equation*}
 \omega(z,A,D):= 0,\qquad z\in A^{\ast}.
\end{equation*}

Geometric  properties  of the harmonic measure will be discussed in   Section
4 below.
By Theorem \ref{thm5.15}  below, if either  $A$  is a Borel set or $D\subset\C,$ then
$\omega(\cdot,A,D)\equiv \omega(\cdot,A^{\ast},D).$

\subsection{Cross and separate holomorphicity}

Let $X,\ Y$  be two complex manifolds of dimension $1,$
  let $D\subset X,$ $ G\subset Y$ be two open sets, let
  $A$ (resp. $B$) be a subset of  $\partial D$ (resp.
  $\partial G$) such that
      $D$ (resp. $G$) is   Jordan-curve-like
   on $A$ (resp. $B$) and that
  $A$ and $B$ are  of positive  length. We define
a {\it $2$-fold cross} $W,$ its {\it regular part} $W^{\ast},$ its {\it  interior} $W^{\text{o}},$ as
\begin{eqnarray*}
W &:=&\X(A,B; D,G)
:=((D\cup A)\times B)\cup (A\times(B\cup G)),\\
W^{\ast} &:=&\X(A^{\ast},B^{\ast}; D,G),\\
W^{\text{o}} &:=&\X^{\text{o}}(A,B; D,G)
:= (A\times  G)\cup (D\times B).
\end{eqnarray*}
Moreover, put
\begin{equation*}
\omega(z,w):=\omega(z,A^{\ast},D)+\omega(w,B^{\ast},G),\qquad
(z,w)\in (D\cup A^{\ast})\times (G\cup B^{\ast}).
\end{equation*}
It  is clear  that $\omega|_{D\times G}$ is harmonic.

For a $2$-fold cross $W :=\X(A,B; D,G)$
define  {\it its wedge}
\begin{equation*}
\widehat{W}=\widehat{\X}(A,B;D,G)
:=\left\lbrace (z,w)\in (D\cup A^{\ast})\times (G\cup B^{\ast}):\ \omega(z,w)  <1
\right\rbrace.
\end{equation*}
Then the set of all interior points of the wedge
$\widehat{W}$ is given by
\begin{equation*}
 \widehat{W}^{\text{o}} :=\widehat{\X}^{\text{o}}(A,B;D,G)
 :=\left\lbrace (z,w)\in D\times G :\  \omega(z,w)<1
\right\rbrace.
\end{equation*}
In particular, if $A$ (resp. $B$) is an open set of $\partial D$  (resp.  $\partial G$), one has
$A\times B\subset A^{*}\times B^{\ast}$ and $W\subset W^{\ast}\subset \widehat{W}.$

We say that a function $f:W\longrightarrow \C$ is {\it separately holomorphic}
{\it on $W^{\text{o}} $} and write $f\in\mathcal{O}_s(W^{\text{o}}),$   if
 for any $a\in A $ (resp.  $b\in B$)
 the function $f(a,\cdot)|_{G}$  (resp.  $f(\cdot,b)|_{D}$ )  is holomorphic  on $G$  (resp. on $D$).

 We say that a function $f:\ W\longrightarrow \C$
 (resp. $f:A\times B\longrightarrow \C$)  is  {\it separately continuous}
{\it on $W $}  (resp. {\it on $A\times B$}) and
write $f\in\mathcal{C}_s(W)$ (resp.  $f\in\mathcal{C}_s(A\times B)$),   if
it is continuous with respect to any variable when the remaining variable is fixed.

In the remaining part of this subsection we introduce two notions.
As in Subsection 2.1  we like to point out that  these notions are intrinsic, i.e., they do not depend
on any choice   we made in their definitions.

 We say that a function $f:\ A\times B\longrightarrow \C$
   is  {\it Jordan-measurable  on $A\times B$},   if  for every point $\zeta\in A$ with type $n$
   (resp.  $\eta\in B$ with type $m$) there is an open neighborhood $U$ of $\zeta$
(resp.  $V$ of $\eta$) such that         $U\cap D=\bigcup\limits_{1\leq j\leq n} U_j$
(resp.    $V\cap G=\bigcup\limits_{1\leq k\leq m} V_k$) with Jordan  domains $U_j, \ V_k,$ and  conformal mappings
$\Phi_j$  (resp.  $\Psi_k$)  from  $U_j$ (resp.  $V_k$) onto $E$ which extends homeomorphically from
$\overline{U_j}$ (resp. $\overline{V_k}$) onto $\overline{E}$ such that
 $f(\Phi_j^{-1}(\cdot),\Psi_k^{-1}(\cdot)):\  \Phi_j
(\overline{U_j}\cap A)\times  \Psi_k
(\overline{V_k}\cap B)\longrightarrow\C $ is  Lebesgue measurable.

Two Jordan-measurable functions  $f,g:\ A\times B\longrightarrow \C$  are said to be  {\it equal a.e. on $A\times B,$}
 if  for every point $\zeta\in A$ with type $n$
   (resp.  $\eta\in B$ with type $m$),    the functions
   \begin{equation*}
   f(\Phi_j^{-1}(\cdot),\Psi_k^{-1}(\cdot)),
   g(\Phi_j^{-1}(\cdot),\Psi_k^{-1}(\cdot)):\  \Phi_j
(\overline{U_j}\cap A)\times  \Psi_k
(\overline{V_k}\cap B)\longrightarrow\C
\end{equation*}
  are equal  a.e. (we keep the previous notation).

 We say that a function $f:\widehat{W}^{\text{o}}\longrightarrow \C$ admits an {\it angular limit}
  $\lambda\in\C$ at $(a,b)\in \widehat{W}$ if the following limit holds:\\
  {\bf  Case 1: $a\in D$ and $b\in G:$}
  \begin{equation*}
\lim\limits_{z\to a, w\to b} f(z,w)=\lambda;
\end{equation*}
{\bf Case 2:  $a\in\ A^{\ast}$  and  $b\in G:$}
 \begin{equation*}
\lim\limits_{z\to a,\  z \in  \mathcal{A}_{\alpha}(a),\ w\to b} f(z,w)=\lambda,\qquad 0<\alpha<\frac{\pi}{2};
\end{equation*}
{\bf Case 3:
$a\in D$ and $b\in B^{\ast}:$}
 \begin{equation*}
\lim\limits_{z\to a,\ w\to b,\ w \in  \mathcal{A}_{\alpha}(b)} f(z,w)=\lambda,\qquad 0<\alpha<\frac{\pi}{2};
\end{equation*}
{\bf Case 4:
$a\in A^{\ast}$ and $b\in B^{\ast}:$}
 \begin{equation*}
\lim\limits_{z\to a,\  z \in  \mathcal{A}_{\alpha}(a),\ w\to b,\ w \in  \mathcal{A}_{\alpha}(b)} f(z,w)=\lambda,
\qquad 0<\alpha<\frac{\pi}{2}.
\end{equation*}
\subsection{Motivations for our work}
We are now able to formulate what, in the sequel,  we  quote  as the {\it
classical  version of the
boundary cross theorem.}

\renewcommand{\thethmspec}{Theorem 1}
\begin{thmspec} (Gonchar \cite{go1,go2})
Let  $D,\ G\subset \C$ be  Jordan domains and  $A $ (resp. $B$) a nonempty open set
of the boundary $\partial D$  (resp. $\partial G$).
Then, for any  function $f\in\mathcal{C}(W)\cap\mathcal{O}_s(W^{\text{o}}),$ there is a unique function
$\hat{f}\in\mathcal{C}(\widehat{W})\cap\mathcal{O}(\widehat{W}^{\text{o}})$
such that $\hat{f}=f$ on $W.$ Moreover, if $\vert f\vert_W<\infty$ then
\begin{equation*}
 \vert \hat{f}(z,w)\vert\leq \vert f\vert_{A\times B}^{1-\omega(z,w)} \vert
 f\vert_W^{\omega(z,w)},\qquad (z,w)\in\widehat{W},
\end{equation*}
where $W$, $W^{\text{o}}$, and $\widehat {W}$ denote the $2$-fold cross, its
interior and its wedge, respectively, associated to  $A,\ B$, $D,\ G$.
\end{thmspec}

Theorem 1 admits various generalizations. The following  result
is announced by Gonchar in \cite{go1}.

\renewcommand{\thethmspec}{Theorem 2}
  \begin{thmspec}
Let  $D,\ G\subset \C$ be  Jordan domains and  $A $ (resp. $B$) a nonempty open and rectifiable set
of the boundary $\partial D$  (resp. $\partial G$).
Let  $f$ be a function defined on the $2$-fold cross   $W$  with the following properties:
\begin{itemize}
\item[(i)] $f|_{W^{\text{o}}}\in\mathcal{C}(W^{\text{o}})\cap \mathcal{O}_s(W^{\text{o}});$
\item[ (ii)]  $f$ is locally bounded on $W;$
\item[ (iii)] for any $a\in  A$  (resp. $b\in B$),
 the holomorphic function $f(a,\cdot)|_{G}$  (resp.  $f(\cdot,b)|_{D}$) has the angular  limit
  $f_1(a,b)$ at $b$ for a.e. $b\in B$  (resp.  $f_2(a,b)$ at  $a$ for a.e. $a\in A$)
and $f_1=f_2=f$  a.e. on  $ A\times B.$
\end{itemize}
1) Then  there is a unique function
$\hat{f}\in\mathcal{O}(\widehat{W}^{\text{o}})$
such that
\begin{equation*}
\lim\limits_{(z,w)\in \widehat{W}^{\text{o}},\ (z,w)\to (\zeta,\eta)}\hat{f}(z,w)=f(\zeta,\eta),\qquad(\zeta,\eta)\in W^{\text{o}}.
\end{equation*}
2) If, moreover, $\vert f\vert_W<\infty,$ then
\begin{equation*}
 \vert \hat{f}(z,w)\vert\leq \vert f\vert_{A\times B}^{1-\omega(z,w)} \vert
 f\vert_W^{\omega(z,w)},\qquad (z,w)\in\widehat{W}^{\text{o}}.
\end{equation*}
 3)  If, moreover,   $f$ is continuous at a
point $(a,b)\in A\times B,$ then
\begin{equation*}
\lim\limits_{(z,w)\in \widehat{W}^{\text{o}},\ (z,w)\to (a,b)}\hat{f}(z,w)=f(a,b) .
\end{equation*}
\end{thmspec}

On the other hand,  the following  result
 due to  Dru\.{z}kowski \cite{dr} gives a different flavor.

\renewcommand{\thethmspec}{Theorem 3}
  \begin{thmspec}
  Let  $D,\ G\subset \C$ be  Jordan domains and  $A $ (resp. $B$) a nonempty open connected set
of the boundary $\partial D$  (resp. $\partial G$).
Let  $f$ be a function defined on  $W$  with the following properties:
\begin{itemize}
\item[(i)] $f\in\mathcal{C}_s(W)\cap \mathcal{O}_s(W^{\text{o}});$
\item[ (ii)]  $f$ is locally bounded on $W;$
\item[(iii)]  $f|_{A\times B}$
is continuous on $A\times B.$
\end{itemize}
Then all conclusions of Theorem 1 still hold.
\end{thmspec}

Observe that all these theorems require  the following very strong hypothesis:
$D$  and $G$  are Jordan domains in $\C$  and $A\times B$ is an  open set of $\partial D\times\partial G.$
 Moreover, the assumptions on the boundedness and continuity of $f$
 are rather restrictive.

A natural question is whether Theorems 1--3 are still true if
$D,\ G$ are    open sets in complex manifolds of dimension $1$  and
the  $A$  (resp.  $B$) is not necessarily an open set of $\partial D$  (resp.  $\partial G$). In addition, if
 one  drops the hypothesis on the
local boundedness and the  continuity of $f,$
 can one obtain a holomorphic  extension   of $f$ and what are its properties?
These matters seem to be of interest, especially when one seeks to
generalize Theorems 1--3 to higher dimensions.

The present paper is motivated by these questions. Our first
purpose  is to generalize Gonchar's theorems to a very general
situation, where $D,\  G$ are, in some sense, almost general open
subsets of complex manifolds of dimension $1$   and where the
boundary sets $A,\  B$ are almost general subsets of $\partial D,\
\partial G.$ Our second goal is to establish, in this general
context, an extension theorem analogous to Dru\.{z}kowski's
theorem with a minimum of hypotheses on $f.$
\section{Statement of the  main results and outline of the proofs}
We are now ready to state our  main result.

\renewcommand{\thethmspec}{Theorem A}
  \begin{thmspec}
  Let $X,\ Y$  be two complex manifolds of dimension $1,$
  let $D\subset X,$ $ G\subset Y$ be two open sets and
  $A$ (resp. $B$)  a subset of  $\partial D$ (resp.
  $\partial G$) such that
      $D$ (resp. $G$) is   Jordan-curve-like
   on $A$ (resp. $B$) and that
  $A$ and $B$ are  of positive  length.
Let  $f:\ W\longrightarrow \C$ be   such that:
\begin{itemize}
\item[ (i)] $f$ is locally bounded on $W$ and  $f\in \mathcal{O}_s(W^{\text{o}});$
\item[(ii)] $f|_{A\times B}$ is Jordan-measurable;
\item[ (iii)]   for any $a\in  A$  (resp. $b\in B$),
 the holomorphic function $f(a,\cdot)|_{G}$  (resp.  $f(\cdot,b)|_{D}$) has the angular  limit
  $f_1(a,b)$ at $b$ for a.e. $b\in B$  (resp.  $f_2(a,b)$ at  $a$ for a.e. $a\in A$)
and   $f_1=f_2=f$  a.e. on $ A\times B.$
\end{itemize}

Then  there exists a unique function
$\hat{f}\in\mathcal{O}(\widehat{W}^{\text{o}})$
with the following property:\\
1)  there are  subsets $\tilde{A}\subset A\cap A^{\ast}$ and $
\tilde{B}\subset B\cap B^{\ast}$  such that
 the sets  $A\setminus  \tilde{A}$ and    $B\setminus  \tilde{B}$
are of zero length\footnote{ Under this condition it follows from   Part 1) of
Theorem \ref{thm5.15} below that  $\tilde{A}\subset
\tilde{A}^{\ast}$  and  $\tilde{B}\subset\tilde{B}^{\ast}$   .}
and
  $\hat{f}$  admits the angular limit $f(\zeta,\eta)$ at every point
  $(\zeta,\eta)\in\X^{\text{o}}(\tilde{A},\tilde{B};
  D,G).$

In addition, $\hat{f} $ enjoys the following properties:\\
2) If $\vert f\vert_W<\infty,$ then
\begin{equation*}
 \vert \hat{f}(z,w)\vert\leq \vert f\vert_{A\times B}^{1-\omega(z,w)} \vert
 f\vert_W^{\omega(z,w)},\qquad (z,w)\in\widehat{W}^{\text{o}}.
\end{equation*}
3) For any $(a_0,w_0)\in A^{\ast}\times G$ (resp. $(z_0,b_0)\in D\times B^{\ast}$)
 if  $\lim\limits_{(z,w)\to (a_0,w_0),\ (z,w)\in W} f(z,w) (=:\lambda)$
 (resp.  $\lim\limits_{(z,w)\to (z_0,b_0),\ (z,w)\in W} f(z,w) (=:\lambda)$)  exists, then $\hat{f}$ admits
 the angular limit $\lambda$ at $(a_0,w_0)$  (resp.  at  $(z_0,b_0)$).
\\
4) For any $(a_0,b_0)\in A^{\ast}\times B^{\ast},$ if  $\lim\limits_{(a,b)\to (a_0,b_0)\ (a,b)\in A\times B} f(a,b)  (=:\lambda)$ exists,
then  $\hat{f}$ admits the angular limit $\lambda  $ at $(a_0,b_0).$
\\
5)
If $f|_{A\times B}$ can be extended to a continuous function  defined
on   $A^{\ast}\times B^{\ast},$ then $f$ can be extended to   a unique  continuous function (still denoted
by) $f$ defined on $W^{\ast}:=\X(A^\ast,B^\ast;D,G)$ and
$\hat{f}$ admits the angular limit $f(\zeta,\eta)$ at every $(\zeta,\eta)\in W^{\ast}$ and
$f_1=f_2=f$ on $(A\cap A^{\ast})\times (B\cap B^{\ast}).$
\end{thmspec}

Theorem A has an immediate consequence.

\renewcommand{\thethmspec}{Corollary A'}
  \begin{thmspec}
We keep the hypotheses and the notation of Theorem A.
Suppose in addition that $f\in\mathcal{C}(W^{\text{o}}).$
Then the  function
$\hat{f}\in\mathcal{O}(\widehat{W}^{\text{o}})$
 provided by  Theorem A admits the angular limit $f(\zeta,\eta)$ at every point
   $(\zeta,\eta)\in ((A\cap A^{\ast})\times G)\cup  (D\times (B\cap B^{\ast})).$
\end{thmspec}

It is worthy to note that  Theorem A and Corollary A' generalize, in some sense, Theorems 1--3.

Now we drop   the  hypothesis
on local boundedness and continuity of $f.$ Then
the examples of Dru\.{z}kowski in \cite{dr} (see  Section 10  below) show that,  without these conditions,
the extended function $\hat{f}$ (if it does exist) is, in general, not continuous
on $\widehat{W}.$ However, our second main result gives a partially
positive answer to this question.

\renewcommand{\thethmspec}{Theorem B}
  \begin{thmspec}
   Let $X,\ Y$  be two complex manifolds of dimension $1,$
     let $D\subset X,$ $ G\subset Y$ be two open sets, let
  $A$ (resp. $B$) be a subset of  $\partial D$ (resp.
  $\partial G$) such that
      $D$ (resp. $G$) is   Jordan-curve-like
   on $A$ (resp. $B$) and that
  $A$ and $B$ are  of positive  length.
Let  $f:\ W\longrightarrow \C$  satisfy   the following properties:

\begin{itemize}
\item[ (i)] $f|_{A\times B}\in \mathcal{C}_s(A\times B)$ and $f\in \mathcal{O}_s(W^{\text{o}});$
\item[(ii)]
 for any $a\in  A$  (resp. $b\in B$),  the function $f(a,\cdot)$  (resp.  $f(\cdot,b)$)
is locally bounded on $G\cup B$ (resp. $D\cup A$)
and
 the (holomorphic)  restriction function $f(a,\cdot)|_{G}$  (resp.  $f(\cdot,b)|_{D}$) has the angular  limit
  $f(a,b)$ at $b$ for every  $b\in B$  (resp.   at  $a$ for every  $a\in A$).
\end{itemize}
Then  there are  subsets $\tilde{A}\subset A\cap A^{\ast}$ and  $
\tilde{B}\subset B\cap B^{\ast},$ and a unique function
$\hat{f}\in\mathcal{O}(\widehat{W}^{\text{o}})$
 with  the following properties:
\\
1)  the sets  $A\setminus  \tilde{A}$ and    $B\setminus  \tilde{B}$
are of zero length;
\\
2) $\hat{f}$  admits the angular limit $f(\zeta,\eta)$ at every point
   $(\zeta,\eta)\in \X(\tilde{A},\tilde{B};D,G).$
\end{thmspec}

Observe that if $f\in \mathcal{C}_s(W)\cap \mathcal{O}_s(W^{\text{o}}),$ then
conditions (i)--(ii) above are fulfilled.
Although our results have been stated only  for the case  of a $2$-fold cross,
they can be formulated   for the general case of an $N$-fold cross  with $N\geq 2$  (see also \cite{nv,pn}).

\smallskip

Now we present some ideas how to prove Theorems A and B.

\smallskip

 Our method consists of two steps.
In the first step we  suppose that  $D$ and $G$ are    Jordan domains in $\C.$
In the second one we treat the general case. The key technique here is to
use {\it level sets} of the harmonic measure. More precisely,
we exhaust $D$  (resp.  $G$) by the  level sets of the harmonic measure
$\omega(\cdot,A,D)$  (resp. $\omega(\cdot,B,G)$),  i.e. by
$D_{\delta}:=\left\lbrace z\in D:\ \omega(z,A,D)<1-\delta \right\rbrace$
(resp.  $G_{\delta}:=\left\lbrace w\in G:\ \omega(w,B,G)<1-\delta \right\rbrace              $)  for  $0<\delta<1.$

In order to carry out the first step,
we improve Gonchar's method \cite{go1,go2} and make intensive  use of  Carleman's formula and of
geometric properties of the level sets of harmonic measures.

The main ingredient  for the second step is a   mixed cross type theorem (see also \cite{pn})
valid for measurable boundary sets in the context of complex manifolds of dimension $1.$ 
We prove that theorem using a recent work of Zeriahi  (see \cite{zr}) and
the classical  method of doubly orthogonal bases of Bergman type.

In the second step we  apply this  mixed cross type theorems  in order  to
prove Theorems A and B  with $D$  (resp. $G$) replaced by $D_{\delta}$  (resp.  $G_{\delta}$). Then we construct
the solution  for the  original open sets $D$ and $G$ by means of a gluing procedure (see also \cite{nv}).

\section{Properties of  the harmonic measure and its level sets}

In this section  $X$ is  a complex manifold of dimension $1,$
   $D\subset X$  an open set,  and
$A$  a nonempty Jordan-measurable subset of $\partial D.$
Observe that then $\partial D$ is non-polar.

Let $\mathcal{P}_D$ be the generalized Poisson  integral of $D.$
If, in addition, $A $ is a Borel set, then, by Theorem 4.3.3 of
\cite{ra},
the harmonic measure of $\partial D\setminus A$ is given by
\begin{equation} \label{eq5.1}
 \omega(\cdot,A,D) =\mathcal{P}_D[1_{\partial D\setminus
 A}].
\end{equation}

The following elementary lemma will be  useful.

\begin{lem}\label{lem5.9}
Let $E$ be the unit disc and  $A$ a measurable subset of
$\partial E.$\\
1)
 Let
$u$ be a subharmonic function defined on $E$ with $u\leq 1$ and let $\alpha\in (0,\frac{\pi}{2})$
be such that
\begin{equation*}
 \limsup_{z\to \zeta,\ z\in\mathcal{A}_{\alpha}(\zeta)}
 u(z)\leq 0\qquad\text{for a.e.}\ \zeta\in A.
\end{equation*}
Then $u\leq \omega(\cdot,A,E)$ on $E.$
\\
2) For all density points $\zeta$ of $A,$
\begin{equation*}
\lim\limits_{z\to \zeta,\ z\in
\mathcal{A}_{\alpha}(\zeta) }\omega(z,A,E)=0,\qquad 0<\alpha<\frac{\pi}{2}.
\end{equation*}
In particular,  all density points of $A$ are contained in $A^{\ast}.$\\
3) For all interior points $\zeta$ of $A,$
\begin{equation*}
\lim\limits_{z\to \zeta }\omega(z,A,E)=0.
\end{equation*}
\end{lem}
\begin{proof}  It follows almost immediately  from the explicit formula  for $\mathcal{P}_E.$
\end{proof}

\begin{prop} (Maximum Principle) Let $u\in\mathcal{SH}(D)$ be such that $u$ is bounded
from the above and
\begin{eqnarray*}
\limsup\limits_{z\to\zeta} u(z)&\leq &   0,\qquad \zeta\in \partial D\setminus A,\\
 \limsup\limits_{z\to \zeta,\ z\in
\mathcal{A}_{\alpha}(\zeta)} u(z)&\leq &   0,\qquad \zeta\in   A,\  0<\alpha<\frac{\pi}{2}.
\end{eqnarray*}
Then $u\leq 0$ on $D.$
\end{prop}
\begin{proof}  Suppose that  $u<M$ for some  $M.$
Let $\zeta_0$ be an arbitrary point of $A.$ Fix a  Jordan domain $U$ such
that $U\subset D$ and
$\partial U\cap\partial D$ is a closed arc which is a neighborhood of $\zeta_0$ in $\partial D.$
Let $B$ be an open arc in $\partial U\cap\partial D$ which contains $\zeta_0.$
Part 1) and Part 3) of   Lemma \ref{lem5.9} applied  to $u|_U$ yield that
\begin{equation*}
 \limsup\limits_{z\to \zeta,\ z\in U
 } u(z)\leq  M\cdot \limsup\limits_{z\to \zeta,\ z\in U
 } \omega(z,B,U)=   0,\qquad \zeta\in B.
\end{equation*}
 Since $\zeta_0\in B$  and $\zeta_0$ is an arbitrary point of $A$,
 we deduce that
 \begin{equation*}
 \limsup\limits_{z\to \zeta,\ z\in
D
 } u(z)\leq      0,\qquad \zeta\in A.
\end{equation*}
 Combining this with the hypothesis,
 the desired conclusion follows from the
 classical Maximum Principle (see Theorem 2.3.1 in \cite{ra}).
\end{proof}

 In the sequel we formulate some important stability property
of  the  harmonic measure.
Let $\phi:\
\partial D\longrightarrow \R$ be a bounded function.
The associated {\it Perron function}  $H_{D,A}:\ D\longrightarrow\R$  is defined by
\begin{equation}\label{Perronfor}
 H_{D,A}[\phi]:=\sup\limits_{u\in \widehat{\mathcal{U}}}u,
\end{equation}
where $\widehat{\mathcal{U}}=\widehat{\mathcal{U}}(\phi,A,D)$ denotes the family of all
subharmonic functions $u$ on
$D$ such that
\begin{eqnarray*}
\limsup\limits_{z\to\zeta} u(z)&\leq &  \phi(\zeta),\qquad \zeta\in \partial D\setminus A,\\
 \limsup\limits_{z\to \zeta,\ z\in
\mathcal{A}_{\alpha}(\zeta)} u(z)&\leq &  \phi(\zeta),\qquad \zeta\in   A,\  0<\alpha<\frac{\pi}{2}.
\end{eqnarray*}
In the sequel,  $\widehat{\mathcal{U}}(A,D)$ will stand for
 $\widehat{\mathcal{U}}( 1_{\partial D\setminus A},A,D).$

Using the above proposition, the corresponding results in   Sections 4.1 and 4.2
of \cite{ra} with respect to $H_{D,A}$ (instead of $H_D$) are still valid making
the obviously  necessary changes. In particular, we have the following
(see Corollary 4.2.6 in \cite{ra}):
\begin{prop}
Let  $D$ be an open subset of $X,$    $A$   a nonempty Jordan-measurable subset of $\partial D,$
 and
  $\phi:\ \partial D\longrightarrow \R$   a bounded function which is
  continuous nearly everywhere\footnote{
  A property is said to hold nearly everywhere on $\partial D$ if
  it holds everywhere on $\partial D\setminus\mathcal{N}$ for some
  Borel polar set $\mathcal{N}.$} on $\partial D.$
Then there exists a unique bounded harmonic function $h$ on $D$ such that
$\lim\limits_{z\to\zeta} h(z)=\phi(\zeta)$
for nearly everywhere $\zeta\in\partial D.$ Moreover,  $h=H_{D}[ \phi]=H_{D,A}[\phi].$
\end{prop}

In virtue of this result, Theorem 4.3.3 in \cite{ra} is still valid in the context
of $H_{D,A}.$ More precisely,
\begin{prop}\label{prop5.13}
Let  $D$ be an open subset of $X,$    $A$   a nonempty Jordan-measurable subset of $\partial D,$
and
$\phi:\ \partial D\longrightarrow \R$ a bounded Borel function.
 Then  $H_D[\phi]=H_{D,A}[\phi]=\mathcal{P}_D[\phi].$
\end{prop}
In the special case  $X:=\C$  we can say  even  more.
\begin{prop}\label{prop5.14}
Let $D$ be a proper  open subset of $\C.$
Let $A$ be a  nonempty Borel subset of $\partial D$ such that  $D$ is Jordan-curve-like on $A$
and  $A$ is of zero length.
Then  $\mathcal{P}_D[1_A]\equiv 0$ on $D.$
 \end{prop}
\begin{proof}
Suppose without loss of generality that  $D$ is Jordan-curve-like on the interval
$[0,1]\subset \partial D$ and that
$A$ is a Borel subset of $[0,1]$ with  $\mes(A)=0.$  Since  $D\subset  \C\setminus [0,1],$
it follows from the Subordination Principle   that
\begin{equation*}
\mathcal{P}_D [1_A]\leq  \mathcal{P}_{\C\setminus [0,1]}[ 1_A]\qquad  \text{on}\ D.
\end{equation*}
Therefore, it suffices to show that $\mathcal{P}_{\C\setminus [0,1]} [1_A] \equiv 0$
on $\C\setminus [0,1].$
To this end consider the conformal mapping $\Phi(z):=\sqrt{\frac{1}{z}-1}$ which
maps $\C\cup\{\infty\}\setminus [0,1]$  onto $\H:=\left\{z\in\C:\  \Im z>0\right\}.$  It
is not difficult to show that
\begin{equation*}
 \mathcal{P}_{\C\setminus [0,1]} [1_A]= \mathcal{P}_{\H} [1_{\Phi(A)} ]\circ \Phi^{-1} \equiv 0.
 \end{equation*}
 This concludes the proof.
\end{proof}
Now  we arrive at one of the main results of the section
\begin{thm}\label{thm5.15}
Let  $D$ be an open subset of $X,$     $A$   a nonempty Jordan-measurable subset of $\partial D,$
and $\mathcal{N}$  a Jordan-measurable subset of $\partial D$  which is of zero length.  \\
1) Then  $A^{\ast}$ is a Borel set and $(A^{\ast})^{\ast}=A^{\ast}$
and $(A\setminus \mathcal{N})^{\ast}=A^{\ast}$  and  $A\setminus A^{\ast}$ is of zero length.\\
2) If $A$ is a Borel set then  $\omega(z,A,D)= H_{D,A} [1_{\partial D\setminus A}]$ for $z\in D.$
In particular,
\begin{equation*}
\omega(z,A^{\ast},D)= H_{D,A^{\ast}} [1_{\partial D\setminus A^{\ast}}]=
 H_{D,(A\cap A^{\ast})\setminus\mathcal{N}}\left[ 1_{\partial D\setminus[(A\cap A^{\ast})\setminus\mathcal{N}]}
 \right]=\omega(z,A,D) ,\qquad z\in D.
 \end{equation*}
3) If  $X=\C$ then  $\omega(z,A,D)=\omega(z,A\setminus\mathcal{N},D)=\omega(z,A^{\ast},D).$
\end{thm}
\begin{proof}
Part  1)  can be checked using the definition  and Lemma \ref{lem5.9}.

Part 2) is an immediate consequence of Proposition \ref{prop5.13} and Part 1).

Now  we turn to Part 3).
 Choose two Borel sets $A_1, \ A_2$ so that $A_1\subset A\setminus\mathcal{N}$ and $ A\subset A_2\subset\partial D$ and
   $A_2\setminus A_1$  is of zero length. Then
  we conclude by  the Subordination Principle and  Proposition  \ref{prop5.14}  that
\begin{equation*}
\omega (z,A_2,D)\leq  \omega(z,A,D)\leq \omega(z,A\setminus\mathcal{N},D)\leq \omega(z,A_1,D)=\omega (z,A_2,D),\qquad z\in D.
\end{equation*}
This proves the first identity.

Since  $A^{\ast}$ is, by Part 1), a Borel set, Part 2) gives that
\begin{equation*}
\omega(z,A^{\ast},D)=H_{D,A^{\ast}}[ 1_{\partial D\setminus A^{\ast}}].
\end{equation*}
Consequently,  $\omega(z,A,D)\leq  \omega(z,A^{\ast},D),  z\in D.$
On the other hand,   let $B$ be a Borel set  such that  $B\subset A\cap A^{\ast}$  and
$A\setminus B$ is of zero length.
Then
\begin{equation*}
\omega(\cdot,A,D)= \omega(\cdot,B,D)=H_{D,B}[ 1_{\partial D\setminus B}]\geq H_{D,A^{\ast}}[
 1_{\partial D\setminus A^{\ast}}]
=\omega(\cdot,A^{\ast},D)\qquad \text{on}\ D.
\end{equation*}
Combining the above estimates, the proof of the last identity in Part 3)   follows.
  \end{proof}

\begin{prop} \label{prop5.16}
Let  $D$ be an open subset of $X$  and     $A$   a nonempty Jordan-measurable subset of $\partial D.$
Let $\left ( D_k \right )^{\infty}_{k=1}$ be a sequence
of   open subsets $D_k$ of $D$    and
 $\left ( A_k\right )^{\infty}_{k=1}$  a sequence of Jordan-measurable subset of $\partial D$ which are also subsets of $A$  such that
\begin{itemize}
\item[(i)] $D_k\subset D_{k+1}$ and $\bigcup_{k=1}^{\infty} D_k=D;$
\item[(ii)]  $A_k\subset A_{k+1}$ and $A_k\subset \partial D\cap\partial
D_k$ and $D_k$ is Jordan-curve-like on $A_k$
and $\bigcup_{k=1}^{\infty} A_k=A;$
 \item[(iii)]  for any point $\zeta\in A$ there is an open neighborhood
 $V=V_{\zeta}$ of $\zeta$ in $\C$ such that $V\cap D$ =$V\cap D_k$ for
 some $k.$
  \end{itemize}
Then
\begin{equation*}
\omega(z, A^{\ast},  D)=
\lim\limits_{k\to\infty}
\omega(z, A_k^{\ast},  D_k),
  \qquad z\in D. 
\end{equation*}
\end{prop}
\begin{proof} Using the Subordination Principle
it is easy to see that  the sequence $\left( \omega(\cdot,  A_k^{\ast},  D_k)\right)_{k=1}^{\infty}$ is decreasing
and    the following limit
\begin{equation*}
u:=\lim\limits_{k\to\infty}
\omega(\cdot,  A_k^{\ast},  D_k)
\end{equation*}
exists and defines a subharmonic  function in $D.$
By the Subordination Principle again, we have  $u\geq \omega(\cdot,A^{\ast},D).$
Therefore, it remains  to establish the converse inequality.
In virtue of (i)--(iii), we  conclude  that
\begin{equation}\label{eq5.23}
\sup\limits_{0<\alpha<\frac{\pi}{2}}\limsup_{z\to \zeta,\
z\in\mathcal{A}_{\alpha}(\zeta)}u =0,\qquad \zeta\in B ,
\end{equation}
where $B:= \bigcup_{k=1}^{\infty} A_k^{\ast}.$

On the other hand, since  $ (A\cap A^{\ast})\setminus B \subset \bigcup\limits_{k=1}^{\infty}
 (  A_k\setminus   A_k^{\ast}
),$
  Part 1) of Theorem \ref{thm5.15} implies that   $ (A\cap A^{\ast})\setminus B$ is of zero length.
Consequently,  
we deduce from (\ref{eq5.23}) and Part 2) of Theorem \ref{thm5.15} that
  $u(z)\leq \omega(z, A^{\ast},D),$ $ z\in D.$  This completes the proof.
\end{proof}

Next, we introduce a notion which will be relevant for our further study.
\begin{defi}\label{endpoint}
Let $D, G\subset X$ be two open sets such that $G\subset D$  and let $\zeta$  be  a point in $\partial D$  such that
$D$ is Jordan-curve-like at $\zeta.$ Then
the point $\zeta$ is said to be an end-point of
$G$ in $D$ if, for every $0<\alpha<\frac{\pi}{2},$ there is an open
neighborhood $U=U_{\alpha}$ of $\zeta$ such that
$U\cap \mathcal{A}_{\alpha}(\zeta)\subset G.$
The set of all end-points of $G$ in $D$ is denoted by $G^D.$
\end{defi}
It is worthy to remark that the above definition is intrinsic.

The remaining part    of this section is devoted to the study of
level sets of the  harmonic measure. We begin with the following
important properties of these sets.
\begin{thm}\label{thm5.19}
Let $D\subset X$ be an open set and  $A$   a  Jordan-measurable
set of $\partial D$ such that $A$ is of positive length.
  Then, for any $0<\epsilon<1,$ the ``$\epsilon$-level set"
\begin{equation*}
D_{\epsilon}:=\left\lbrace z\in D:\ \omega(z,A^{\ast}, D)<1-\epsilon \right\rbrace
\end{equation*}
enjoys the following properties:
\begin{itemize}
\item[(i)] Let $G_1,$ $G_2$ be  arbitrary distinct connected components of $D_{\epsilon},$
 then $G_1^D\cap G_2^D=\varnothing.$
\item[(ii)] For any point  $\zeta\in A^{\ast},$ there is exactly one connected component
$G$ of $D_{\epsilon}$ such that $\zeta\in G^D.$

\item[(iii)] $G^D\cap A$ is  Jordan-measurable  (on $\partial D$) and of positive length for  every connected
component
$G$ of $D_{\epsilon}.$
\end{itemize}
\end{thm}
\begin{proof}
 To prove  (i), suppose, in order to reach a contradiction, that $  G^D_1\cap G^D_2\not=\varnothing.$
Fix a point $\zeta_0\in G^D_1\cap G^D_2.$ Then, for every $0<\alpha<\frac{\pi}{2},$ there is an
open neighborhood
$ U_{\alpha}$ of $\zeta_0$ such that $\mathcal{A}_{\alpha}(\zeta)\cap U_{\alpha}\subset G_1\cap
G_2.$ This implies that $ G_1\cap
G_2\not=\varnothing.$ Hence, $G_1=G_2,$ which contradicts the hypothesis that
$G_1\not =G_2.$  The proof of  (i) is complete.

 Next, we turn to the proof of  (ii). Fix a  $\zeta_0\in A^{\ast}.$ In virtue of assertion
 (i), it suffices to show the existence of a connected component $G$ of
 $D_{\epsilon}$ such that $\zeta_0\in G^D.$  Since  $\zeta\in A^{\ast},$ for every
 $0<\alpha<\frac{\pi}{2},$ there is an open neighborhood $U_{\alpha}$ of
 $\zeta_0$ such that
 \begin{equation}\label{eq5.24}
\mathcal{A}_{\alpha}(\zeta_0)\cap U_{\alpha}\subset D_{\epsilon}.
 \end{equation}
 Fix an arbitrary $0<\alpha_0<\frac{\pi}{2},$ and let $G$ be the
 connected component of $D_{\epsilon}$ containing $\mathcal{A}_{\alpha_0}(\zeta_0)\cap
 U_{\alpha_0}.$ Since
\begin{equation*}
\Big(\mathcal{A}_{\alpha_0}(\zeta_0)\cap U_{\alpha_0}\Big)\cap
\Big(  \mathcal{A}_{\alpha}(\zeta_0)\cap U_{\alpha}
\Big)\not=\varnothing,\quad 0<\alpha<\frac{\pi}{2},
 \end{equation*}
we deduce from (\ref{eq5.24}) that $G$ also contains $\mathcal{A}_{\alpha}(\zeta_0)\cap
U_{\alpha}$ for every  $0<\alpha<\frac{\pi}{2}.$ Hence $\zeta_0\in G^D.$
The proof of   (ii) is finished.

Finally, we prove   (iii). First, we may find a sequence $(U_k)_{k=1}^{\infty}$ of open
sets of $X$ such
that $U_k\cap D$ is  either a   Jordan domain or the disjoint union of two
Jordan domains and $A\subset
\bigcup\limits_{k=1}^{\infty}\partial(U_k\cap D).$ Since $A$ is
Jordan-measurable, we see that in order to prove the  Jordan-measurability of $G^D\cap A,$  it is
sufficient to
check that $G^D\cap \partial (D\cap U_k)$ is  Jordan-measurable for
every $k\geq 1.$ To prove the latter assertion, fix an $k_0\geq 1$ and let $U:=U_{k_0}.$ Let $\Phi$ be a conformal mapping
  from $D\cap U$ onto $E$ which extends to a homeomorphic mapping (still denoted by) $\Phi$ from
    $\overline{D\cap U}$ onto $\overline{E}.$
It is clear that for any $\zeta\in \partial (D\cap U),$
$\zeta\in G^D$ if and only if $\Phi(\zeta)\in [\Phi(G\cap U)]^E.$
We shall prove, in the sequel, that $[\Phi(G\cap U)]^E$ is a Borel subset of $\partial E.$
 Taking this for granted, then  $G^D\cap \partial (D\cap U)$ is also a
 Borel set. Consequently, $G^D\cap A$ is  Jordan-measurable.

To check that  $[\Phi(G\cap U)]^E$ is a Borel set, put
\begin{equation} \label{eq5.25}
\mathcal{A}_{n,m}(\eta):=\left  \lbrace  w\in E\cap
\mathcal{A}_{\left(1-\frac{1}{n}\right)\cdot\frac{\pi}{2}
}(\eta):\  \vert w- \eta\vert <\frac{1}{m}\right\rbrace,\;  n,m\geq 1,\; \eta\in\partial E.
\end{equation}
For any $n,m,p\geq 1,$ let
\begin{multline}\label{eq5.26}
T_{nmp}:=\Big\lbrace  \eta\in\partial E:\
\mathcal{A}_{n,m}(\eta)\subset \Phi(G\cap U)\quad\text{and}   \\
 \quad \omega(\Phi^{-1}(w),A^{\ast},D)\leq
1-\epsilon-\frac{1}{p}, \forall
w\in \mathcal{A}_{n,m}(\eta) \Big\rbrace.
\end{multline}

We observe the following {\bf geometric fact:}

\smallskip

 {\it Let $\eta_0\in\partial E$ and $(\eta_q)_{q=1}^{\infty}\subset \partial E$
such that $\lim\limits_{q\to\infty}\eta_q=\eta_0.$ Then
\begin{equation*}
 \mathcal{A}_{n,m}(\eta_0)\subset \bigcup\limits_{q=1}^{\infty}
 \mathcal{A}_{n,m}(\eta_q).
\end{equation*}
}

The proof of this fact follows immediately from the geometric shape of the cone $
\mathcal{A}_{n,m}(\eta)$ given in (\ref{eq5.25}).

Let  $(\eta_q)_{q=1}^{\infty}\subset T_{nmp}$  such that  $\lim\limits_{q\to\infty}\eta_q=
\eta_0\in\partial E.$
Using the above geometric fact, we see that $\mathcal{A}_{n,m}(\eta_0)\subset
\Phi(G\cap U).$ This, combined with (\ref{eq5.26}) and the continuity of
$\omega(\Phi^{-1}(\cdot),A,D)|_E,$ implies that $\eta_0\in T_{nmp}.$
Hence, the set
$T_{nmp}$ is closed. Clearly, we have
\begin{equation*}
[\Phi(G\cap U)]^E=\bigcap\limits_{n=1}^{\infty}\bigcup\limits_{m=1}^{\infty}
\bigcup\limits_{p=1}^{\infty}T_{nmp}
\end{equation*}
It follows immediately from this identity that $[\Phi(G\cap U)]^E$ is a
Borel set. Consequently, as was already discussed before, $G^D\cap A$   is  Jordan-measurable.

  To finish assertion (iii), it remains to prove that
  $ G^E\cap A$  is of positive length.
  Suppose, in order to reach a contradiction, that
$ G^E\cap A$  is of zero length.
Consider the following function
\begin{equation*}
 u(z):=
\begin{cases}
\omega(z,A^{\ast},D),
  &z\in   D\setminus G\\
 1-\epsilon, & z \in G
\end{cases}.
\end{equation*}
Then clearly $u\in\mathcal{SH}(D)$   and $u\leq 1.$
In virtue of assertions (i) and (ii) and the definition of locally regular points, we have that
\begin{equation*}
\sup\limits_{0<\alpha<\frac{\pi}{2}}\limsup\limits_{z\to \zeta,\ z\in
\mathcal{A}_{\alpha}(\zeta) }u(z)=\sup\limits_{0<\alpha<\frac{\pi}{2}}\limsup\limits_{z\to \zeta,\ z\in
\mathcal{A}_{\alpha}(\zeta) } \omega(z,A^{\ast},E)=      0,\qquad   \zeta\in (A\cap A^{\ast})\setminus
(G^D\cap A).
\end{equation*}
Consequently, using the notation in (\ref{Perronfor}), we conclude that
$$
u\in \widehat{\mathcal{U}}\left((A\cap A^{\ast})\setminus\mathcal{N},D\right),
$$
 where $\mathcal{N}:=G^D\cap A.$
Since, by our above assumption,
$ \mathcal{N}$  is of zero length, it follows from Theorem
\ref{thm5.15} that   $u\leq \omega(\cdot,A^{\ast},D).$ But on the other hand, one has
$\omega(z,A^{\ast},D)<1-\epsilon=u(z)$ for $z\in G.$ This leads to the desired
contradiction. Hence, the proof of    (iii) is finished.
\end{proof}

\begin{thm}\label{thm5.20}
Let $D\subset X$ be an open set and  $A$   a  Jordan-measurable
set of $\partial D$ such that $A$ is of positive length.
  For any $0\leq \epsilon<1,$ let  $D_{\epsilon}:=\left\lbrace z\in D:\ \omega(z,A^{\ast},D)<1-\epsilon
  \right\rbrace.$\\
  1) For any Jordan-measurable subset $\mathcal{N}\subset \partial D$  of zero length, let
\begin{multline*}
\mathcal{U}_{\epsilon}(A,\mathcal{N},D):=\\
  \left\lbrace u\in\mathcal{SH}(D_{\epsilon}):\ u\leq
1\ \text{and}\ \sup\limits_{0<\alpha<\frac{\pi}{2}}\limsup\limits_{z\to\zeta,\  z\in
\mathcal{A}_{\alpha}(\zeta)}  u(z)\leq  0,\ \zeta\in (A\cap
A^{\ast})\setminus\mathcal{N}
\right\rbrace.
\end{multline*}
Then
$\mathcal{U}_{\epsilon}(A,\mathcal{N},D)=\mathcal{U}_{\epsilon}(A,\varnothing,D).$\\
2) Define the ``harmonic measure of
the $\epsilon$-level set" $\omega_{\epsilon}(\cdot,A,D)$ as
\begin{equation*}
\omega_{\epsilon}(z,A,D) :=\begin{cases}
\sup\limits_{u\in \mathcal{U}_{\epsilon}(A,\varnothing,D)} u(z),
  & z\in    D_{\epsilon} \\
 0, & z \in A^{\ast}
\end{cases}.
\end{equation*}
Then
\begin{equation*}
\omega_{\epsilon}(z,A,D)=\frac{\omega(z,A^{\ast},D)}{1-\epsilon},\qquad z\in D_{\epsilon}\cup A^{\ast}.
\end{equation*}
\end{thm}
\begin{proof}
 Clearly, by definition, $\mathcal{U}_{\epsilon}(A,\varnothing,D)\subset
 \mathcal{U}_{\epsilon}(A,\mathcal{N},D).$  To prove the converse
 inclusion, fix an arbitrary $u\in \mathcal{U}_{\epsilon}(A,\mathcal{N},D).$
 Consider  the following function
\begin{equation*}
 \hat{u}(z):=
\begin{cases}
\max\left\lbrace (1-\epsilon)u(z),\omega(z,A^{\ast},D)\right\rbrace,
  &z\in   D_{\epsilon}\\
  \omega(z,A^{\ast},D), & z \in D\setminus D_{\epsilon}
\end{cases}.
\end{equation*}
Then   $\hat{u}\in\mathcal{SH}(D)$   and $\hat{u}\leq 1.$ Moreover,
in virtue of   (ii) of Theorem \ref{thm5.19}, we have that
$A^{\ast}\subset (D_{\epsilon})^D.$ Consequently,  for every $\zeta\in
(A\cap A^{\ast})\setminus \mathcal{N},$
\begin{multline}\label{eq5.27}
\sup\limits_{0<\alpha<\frac{\pi}{2}}\limsup\limits_{z\to \zeta,\ z\in
\mathcal{A}_{\alpha}(\zeta) }\hat{u}(z)\\
\leq \max\left\lbrace \sup\limits_{0<\alpha<\frac{\pi}{2}}\limsup\limits_{z\to \zeta,\ z\in
\mathcal{A}_{\alpha}(\zeta) } u(z),
\sup\limits_{0<\alpha<\frac{\pi}{2}}\limsup\limits_{z\to \zeta,\ z\in
\mathcal{A}_{\alpha}(\zeta) }  \omega(z,A,D) \right\rbrace.
\end{multline}
Observe that  the first term in the latter line of (\ref{eq5.27}) is equal
to $0$ because $u\in\mathcal{U}_{\epsilon}(A,\mathcal{N},D).$ In addition,
  the second term in the latter line of (\ref{eq5.27}) is
 also equal to $0.$ Hence, $\hat{u}\in\widehat{\mathcal{U}}\left( (A\cap A^{\ast})\setminus \mathcal{N},D\right).$
Consequently, by Theorem  \ref{thm5.15},  $\hat{u}\leq\omega(\cdot,A^{\ast},D).$
In particular, one  has
\begin{equation}\label{eq5.28}
u(z)\leq \frac{\omega(z,A^{\ast},D)}{1-\epsilon},\qquad z\in D,\ u\in\mathcal{U}_{\epsilon}(A,\mathcal{N},D).
\end{equation}
On the other hand, it is clear that
$\frac{\omega(\cdot,A^{\ast},D)}{1-\epsilon} \in\mathcal{U}_{\epsilon}(A,\varnothing,D)\subset
\mathcal{U}_{\epsilon}(A,\mathcal{N},D).$
This, combined with (\ref{eq5.28}), implies the desired conclusions of Part 1) and Part 2).
\end{proof}
An immediate consequence of Theorem \ref{thm5.20} is the following
Two-Constant Theorem for level sets.
\begin{cor}\label{cor5.21}
Let $D\subset X$ be an open set and  $A,\  \mathcal{N}$   two  Jordan-measurable
subsets of $\partial D$ such that $A$ is of positive length and $\mathcal{N}$  is of zero length.
Let $0\leq \epsilon<1$ and put  $D_{\epsilon}:=\left\lbrace z\in D:\ \omega(z,A^{\ast},D)<1-\epsilon
  \right\rbrace.$
    If $u\in\mathcal{SH}(D_{\epsilon})$ satisfies $u\leq M$ on
    $D_{\epsilon}$ and $\sup\limits_{0<\alpha<\frac{\pi}{2}}\limsup\limits_{z\to \zeta,\ z\in
\mathcal{A}_{\alpha}(\zeta) } u(z)\leq m,$  $\zeta\in( A\cap
A^{\ast})\setminus \mathcal{N},$ then
\begin{equation*}
u(z)\leq m(1-\omega_{\epsilon}(z,A,D))+M\cdot\omega_{\epsilon}(z,A,D).
\end{equation*}
\end{cor}
\section{Boundary behaviour of the Gonchar--Carleman operator}
Before recalling the   Gonchar--Carleman
operator and investigating its boundary behavior, we first  introduce  the following notion and study its
properties.
\subsection{Angular Jordan domains}
Let $E$ be the unit disc. We begin with the
\begin{defi} \label{def6.1}
For  every closed subset $F$ of $\partial E$ and any real number $h$ such that
  $\mes(F)>0$
and $\sup_{x,y\in F}\vert x-y\vert<h<1-\frac{\sqrt{2}}{2},$ the open set
\begin{equation*}
\Omega=\Omega(F,h):=\bigcup\limits_{\zeta\in F}
\left\lbrace z\in \mathcal{A}_{\frac{\pi}{4}}(\zeta):\ \vert z\vert>1-h  \right\rbrace
\end{equation*}
is called the {\it angular Jordan domain} with {\it base} $F$ and {\it height} $h.$
\end{defi}
Now we give   a list of properties of such  angular Jordan domains.
\begin{prop} \label{prop6.2}
Let $\Omega=\Omega(F,h)$ be an angular Jordan domain.\\
1) Then there exist exactly two points $\zeta_1,\ \zeta_2\in F$
such that $\vert \zeta_1-\zeta_2\vert=\sup_{x,y\in F}\vert x-y\vert$ and
  $F\subset [\zeta_1,\zeta_2],$ where $[\zeta_1,\zeta_2]$ is the (small) closed arc of $\partial E$
  which is oriented in the positive sense and which starts from
 $\zeta_1$ and ends at  $\zeta_2.$ \\
2) Write the open set  $ [\zeta_1,\zeta_2]\setminus F$  as the union of disjoint open arcs
\begin{equation*}
[\zeta_1,\zeta_2]\setminus F=\bigcup_{j\in J}  (a_j,b_j),
\end{equation*}
where $(a_j,b_j)$ is the (small) open arc of $\partial E$ which goes from
  $a_j$  to   $b_j$ and which is oriented in the positive sense, and the index set $J$ is finite
  or countable.

For $j\in J,$ we construct the isosceles triangle with the three vertices $a_j,$ $b_j$ and
$c_j$ such that the base of the isosceles triangle    is the segment
connecting $a_j$ to $b_j,$ and $c_j$ satisfies
\begin{equation*}
 \arg\left( \frac{c_j-a_j}{a_j} \right)=\frac{3\pi}{4}\quad\text{and}\quad
 \arg\left( \frac{c_j-b_j}{b_j} \right)=\frac{-3\pi}{4}.
\end{equation*}
Let $[a_j,c_j]$ (resp. $[c_j,b_j]$) denote the segment connecting $a_j$ to
$c_j$ (resp. the segment connecting $c_j$ to $b_j$).
Put
\begin{equation*}
F_0:=F\cup \bigcup_{j\in J} ([a_jc_j]\cup [c_jb_j]).
\end{equation*}
Then $F_0$ is a rectifiable Jordan curve starting from $\zeta_1$ and
ending at $\zeta_2.$\\
3) Let $\eta_1$ (resp. $\eta_2$) be the unique point in the circle
$\partial \B(0,1-h)$ such that
\begin{equation*}
 \arg\left( \frac{\eta_1-\zeta_1}{\zeta_1} \right)=\frac{-3\pi}{4}\quad \Big(\text{resp.}\quad
 \arg\left( \frac{\eta_2-\zeta_2}{\zeta_2} \right)=\frac{3\pi}{4} \Big)
\end{equation*}
and that $\vert \eta_1-\zeta_1\vert$ (resp. $\vert\eta_2- \zeta_2\vert$) is minimal.
Let $ F_1$ (resp. $ F_2$) denote the segment connecting $\eta_1$ to
$\zeta_1$ (resp. the segment connecting $\zeta_2$ to $\eta_2$).
Let $F_3$ be the (small) closed arc of the circle $\partial \B(0,1-h)$ which
starts from $\eta_2$ and ends at $\eta_1$ and which is oriented in the
negative sense.

  Then $\Omega$ is a rectifiable Jordan domain and its boundary
$\Gamma$ consists of the rectifiable Jordan curve $F_0,$ two segments
$F_1,$ $F_2$ and the closed arc $F_3.$\\
4) For every $\epsilon\in (0,\tfrac{h}{4})$  define the dilatation
$\tau_{\epsilon}: \ E\longrightarrow E$  as follows
\begin{equation*}
\tau_{\epsilon}(z):= (1-\epsilon)z,\qquad z\in E.
\end{equation*}
Put
\begin{equation*}
\Omega_{\epsilon}:=\tau_{\epsilon}( \Omega)
\setminus \overline{\B\left(0,(1+\epsilon)(1-h)\right)}.
\end{equation*}
Then $\Omega_{\epsilon}$ is a rectifiable Jordan domain  and its boundary
$\Gamma_{\epsilon}$ consists of
the rectifiable Jordan curve $F_{0\epsilon}:=\tau_{\epsilon}( F_0) ,$
a sub-segment    $F_{1\epsilon}$ of  $\tau_{\epsilon}(F_1),$  a sub-segment   $F_{2\epsilon}$
of $\tau_{\epsilon}(F_2),$
and a closed  arc  $F_{3\epsilon}$ of  $\partial \B\left(0,(1+\epsilon)(1-h)\right).$\\
5) Consider the projection $\tau:\ E\setminus\{0\}\longrightarrow \partial E$ given by
$\tau(z):=\frac{z}{\vert z\vert},$ $z\in  E\setminus\{0\}.$
For every  $\epsilon\in (0,\tfrac{h}{4})$
notice that $F_{0\epsilon}\cup F_{1\epsilon}\cup F_{2\epsilon}=\Gamma_{\epsilon}\setminus \partial
\B\left(0,(1+\epsilon)(1-h)\right).$ Then the two maps
\begin{eqnarray*}
 F_{0\epsilon}\cup F_{1\epsilon}\cup F_{2\epsilon} \ni\zeta & \mapsto &
\tau(\zeta)\in\partial E,\\
 F_{3\epsilon} \ni\zeta & \mapsto  & \tau(\zeta)\in\partial E,
\end{eqnarray*}
are one-to-one. In addition, for any linearly measurable subset $A$ of $\Gamma_{\epsilon},$
\begin{equation*}
\mes(A)\leq 10\cdot \mes(\tau(A)).
\end{equation*}
6) $\Omega_{\epsilon}\nearrow \Omega$ as $\epsilon \searrow 0.$\\
7) For any closed Jordan curve $\mathcal{C}$  contained in  $\Omega$
there is an $\epsilon>0$ such that $\mathcal{C}\subset \Omega_{\epsilon}.$
\\
8) $ \mes(F\setminus \Omega^E)=0.$
\end{prop}
\begin{proof} All assertions are quite simple using an elementary geometric
argument. Therefore, we leave the details of their proofs  to the reader.
However, we will  give   the proof that $\Omega$ is a domain.
This proof will clarify Definition \ref{def6.1}.

In virtue of the condition on $F$ and $h$ given in Definition
\ref{def6.1},  we see that $\left\lbrace z\in \mathcal{A}_{\frac{\pi}{4}}(\zeta):\ \vert z\vert>1-h
\right\rbrace,$  $\zeta\in\partial E,$ is  connected, and that
\begin{multline*}
\left\lbrace z\in \mathcal{A}_{\frac{\pi}{4}}(\zeta):\ \vert z\vert>1-h  \right\rbrace
\cap \left\lbrace z\in \mathcal{A}_{\frac{\pi}{4}}(\eta):\ \vert z\vert>1-h
\right\rbrace\not=\varnothing,\\
\forall \zeta,\eta\in\partial E:\
\vert \zeta-\eta\vert<h<1-\frac{ \sqrt{2}}{2}.
\end{multline*}
Hence, $\Omega$ is a domain.
\end{proof}

\begin{thm}\label{thm6.3}  Let $X$ be a complex manifold of dimension $1,$
 $D\subset X$  an open set and  $A$   a  Jordan-measurable
subset of $\partial D$ such that $A$ is of positive length.
Then, for any $0\leq\epsilon<1$ and  any connected component
$G$ of $D_{\epsilon}:=\left\lbrace z\in D:\  \omega(z,A^{\ast},D)<1-\epsilon \right\rbrace,$
there are an open set $U\subset X,$  a conformal mapping $\Phi:\ E\longrightarrow X$, and an angular Jordan domain $\Omega
=\Omega(F,h)$ such that
\begin{itemize}
\item [(i)] $U\cap D$ is either  a  Jordan domain or the disjoint union of
two  Jordan domains;
\item[(ii)] $\Phi$ maps $E$ conformally onto one  connected component of $ U\cap D$
(notice that, in virtue of (i), $U\cap D$ has at most two connected
components);
\item[(iii)]
 $\Phi(F)\subset A\cap A^{\ast}\cap G^D$ and $\Phi(\Omega)\subset G.$
 \end{itemize}
\end{thm}
\begin{proof}
  We have already shown in the proof of (iii) of Theorem \ref{thm5.19} that there is
   a sequence $(U_k)_{k=1}^{\infty}$ of open sets of $X$ such
that $U_k\cap D$ is  either a   Jordan domain or the disjoint union of two
Jordan domains, and $A\subset
\bigcup\limits_{k=1}^{\infty}\partial(U_k\cap D),$ and $A\cap A^{\ast}\cap G^D$  is of positive length.
 Consequently, there is an index $k_0$ such that
 \begin{equation}\label{eq6.1}
\Big(A\cap A^{\ast}\cap G^D\cap \partial (D\cap U)\Big)\quad\text{is of positive length},
\end{equation}
 where $U:=U_{k_0}.$ Suppose without loss of generality that
 $U\cap D$ is a  Jordan domain. The remaining case where
  $U\cap D$ is the disjoint union of two  Jordan domains may be
 proved in the same way.
Let $\Phi$ be a conformal mapping from $E$ onto $D\cap U.$ By Carath\'eodory Theorem  (see \cite{go}), $\Phi$
extends to a homeomorphic map (still denoted by) $\Phi$ from
$\overline{E}$ onto $\overline{D\cap U}.$  Hence, (i) and (ii) are satisfied.

 On the other hand, it follows from (\ref{eq6.1}) that
\begin{equation}\label{eq6.2}
\mes\Big(\Phi^{-1}\left(A\cap A^{\ast}\cap G^D\cap \partial (D\cap U)\right)\Big)>0.
\end{equation}

For any $m\geq 1,$ let
\begin{equation}\label{eq6.3}
A_m:=\left\lbrace\eta\in\partial E:\ \mathcal{A}_{2,m}(\eta)\subset \Phi^{-1}(G)\right\rbrace,
\end{equation}
where $ \mathcal{A}_{2,m}(\eta)$ is given by formula (\ref{eq5.25}).

Using the Geometric fact just after (\ref{eq5.26}), we see that $A_m$ is
closed. On the other hand, it is clear that
\begin{equation*}
\Phi^{-1}\left(A\cap A^{\ast}\cap G^D\cap \partial (D\cap U)\right)\subset
\bigcup\limits_{m=1}^{\infty}A_m.
\end{equation*}
Therefore, in virtue of (\ref{eq6.2}), there is an index $m_0$ such that
\begin{equation*}
\mes\Big(A_{m_0}\cap \Phi^{-1}\left(A\cap A^{\ast}\cap G^D\cap \partial (D\cap U)\right)\Big)>0.
\end{equation*}
Put $h:=\frac{1}{2m_0}.$ By the latter estimate one may find a
closed set $F$ contained in   $A_{m_0}
\cap \Phi^{-1}\left(A\cap A^{\ast}\cap G^D\cap \partial (D\cap
U)\right)$ such that $\mes(F)>0$ and $\sup\limits_{x,y\in F}\vert x-y\vert
<h.$ Since $h=\frac{1}{2m_0},$ a geometric argument shows that
\begin{equation*}
\left\lbrace z\in \mathcal{A}_{\frac{\pi}{4}}(\zeta):\ \vert z\vert>1-h  \right\rbrace
 \subset  \mathcal{A}_{ 2,m_0}(\zeta) ,\qquad \zeta\in\partial E.
\end{equation*}
This together   with (\ref{eq6.3}) implies that
$\Omega=\Omega(F,h)\subset \Phi^{-1}(G).$ Hence,  (iii) is
verified. This completes the proof.
\end{proof}

In the sequel, the following uniqueness theorem will play a vital role.
\begin{thm}
\label{thm6.4}  Let $X$ be a complex manifold of dimension $1,$
 $D\subset X$  an open set, and  $A,\  \mathcal{N}$   two  Jordan-measurable
subsets of $\partial D$ such that $A$ is of positive length and $\mathcal{N}$  is of zero length.
Let $0\leq\epsilon<1$ and  $G$   a connected component
  of $D_{\epsilon}:=\left\lbrace z\in D:\  \omega(z,A^{\ast},D)<1-\epsilon \right\rbrace.$
If   $f\in\mathcal{O}(G)$   admits the angular
limit $0$ at every point of
 $(A\cap A^{\ast}\cap G^D)\setminus \mathcal{N},$  then $f\equiv 0.$
\end{thm}
\begin{proof} Applying Theorem \ref{thm6.3} we obtain an open set $U$ in
$X,$ a conformal mapping $\Phi$ from $E$ onto $D\cap U$ which extends homeomorphically to
$\overline{E},$  and an angular Jordan domain $\Omega:=\Omega(F,h)$
satisfying assertions (i)--(iii) listed in that theorem.

Consider the function $f\circ \Phi:\ \Omega\longrightarrow\C.$ By the
hypothesis,
$f\circ\Phi\in\mathcal{O}(\Omega)$ admits the angular limit $0$ at a.e
point in $F.$ Since $\mes(F)>0,$  Privalov's Uniqueness  Theorem (see \cite{go}) gives that
$f\circ\Phi\equiv 0$ on $\Omega.$ Hence, $f\equiv 0$ on the subdomain $\Phi(\Omega)
$ of $G.$ This proves  $f\equiv 0.$
\end{proof}
\subsection{Main result of the section}

Let  $D,\ G \subset \C$ be   open  discs
and let  $ A$ (resp.  $B$) be a  measurable subset of $\partial D$  (resp. $\partial G$) with
  $\mes(A)>0$  (resp.  $ \mes(B)>0$). Let  $f$ be a function defined on $W:=\X(A,B;D,G)$  with the following properties:
\begin{itemize}
\item[ (i)]  $f|_{A\times B}$ is measurable and there is a finite constant $C$ with $ \vert f\vert_W<C;$
\item[ (ii)] $f\in \mathcal{O}_s(W^{\text{o}});$
\item[(iii)] there exist two functions $f_1,\  f_2: \ A\times B\longrightarrow \C$
such that for any $a\in A $  (resp.  $b\in B$),
   $f(a,\cdot)$  (resp.  $f(\cdot,b)$) has the angular
 limit $f_1(a,b)$ at $b$ for a.e. $b\in B$  (resp. $f_2(a,b)$ at $a$ for  a.e. $a\in A$),
and  $f_1=f_2= f$  a.e. on $ A\times B.$
\end{itemize}

Let  $\tilde{\omega}(\cdot, A,D)$  (resp.  $\tilde{\omega}(\cdot, B,G)$)  be the conjugate harmonic  function of
$\omega(\cdot,A,D)$ (resp. $\omega(\cdot,B,G)$ ) such that $\tilde{\omega}(z_0,A,D)=0$  (resp.  $\tilde{\omega}(w_0,B,G)=0$)
for a certain fixed point $z_0\in D$  (resp.  $w_0\in G$). Thus we  define the holomorphic functions
$g_1(z):=\omega(z,A,D)+i\tilde{\omega}(z,A,D),$    $g_2(w):=\omega(w,B,G)+i\tilde{\omega}(w,B,G),$         and
\begin{equation*}
g(z,w):=g_1(z)+g_2(w),\qquad (z,w)\in
D\times G.
\end{equation*}

Each function $e^{-g_1}$  (resp.  $e^{-g_2}$) is bounded on $D$  (resp. on  $G$).
 Therefore, in virtue of \cite[p. 439]{go},  we may define  $e^{-g_1(a)}$
  (resp.   $e^{-g_2(b)}$)
for a.e. $a\in A$ (resp. $b\in B$) to be the angular boundary limit of
 $e^{-g_1}$ at $a$   (resp. $e^{-g_2}$  at $b$).

In virtue of (i), for each positive integer $N,$
we define the {\it Gonchar--Carleman  operator} as follows
\begin{equation}\label{eq4.1}
K_N(z,w)=K_N[f](z,w):=\frac{1}{(2\pi
i)^2}\int\limits_{A\times B}e^{-N(g(a,b)-g(z,w))}\frac{f(a,b)da db}{(a-z)(b-w)},\;  (z,w)\in D\times G.
\end{equation}

We recall from Gonchar's work in  \cite{go2}  that  the following limit
\begin{equation}
\label{eq1_Gonchar_Carleman}
K(z,w)=K[f](z,w):=\lim\limits_{N\to\infty}K_N(z,w)
\end{equation}
exists for all $(z,w)\in \widehat{W}^{\text{o}},$ and
its limit is uniform  on compact subsets of
 $\widehat{W}^{\text{o}}.$\\

 The boundary behavior of
Gonchar--Carleman operator is described   below.

\begin{thm}
\label{thm6.5}
We keep the above hypothesis and notation.
 Let  $0< \delta <1,$  $ w\in G$  be  such that
$ \omega(w,B,G)<\delta,$
 and let $U$ be   any connected component  of
 \begin{equation*}
  D_{\delta}:=
 \left\lbrace  z\in D:\ \omega(z,A,D)<1-\delta\right\rbrace.
 \end{equation*}
  Then there
 is an angular Jordan
 domain $\Omega=\Omega(F, h)$
  such that $\Omega\subset U,$  $F\subset A\cap A^{\ast}\cap U^{D},$ and
 the Gonchar--Carleman operator $K[f]$ (see formula (\ref{eq4.1})--(\ref{eq1_Gonchar_Carleman}) above) satisfies
\begin{equation*}
\lim\limits_{z\to a,\ z\in\mathcal{A}_{\alpha}(a)
}K[f](z,w)=f(a,w),\qquad 0< \alpha<\frac{\pi}{2},
\end{equation*}
 for a.e. $a\in F.$
\end{thm}

 The proof of this theorem will be given in Subsection 5.4 below.
\subsection{Preparatory results} For the proof of Theorem
\ref{thm6.5} we need the following results.

In the sequel, for every function $f\in L^1(\partial E,\vert d\zeta\vert),$
 let $ \mathcal{C}[f]$ denote the Cauchy integral
\begin{equation*}
 \mathcal{C}[f](z):=\frac{1}{2\pi i}\int\limits_{\partial E} \frac{f(\zeta)
 d\zeta}{z-\zeta},\qquad z\in E.
\end{equation*}
For a function $F: E\longrightarrow\C,$ the radial maximal function
$M_{\rad} F:\ \partial E\to[0,\infty]$ is defined by
\begin{equation*}
 (M_{\text{rad}}F)(\zeta):=\sup\limits_{0\leq r< 1} \vert F(r\zeta)\vert ,\qquad \zeta\in \partial E.
\end{equation*}
Now we are able to state the following classical result (see Theorem 6.3.1
in Rudin's book \cite{ru})
\begin{thm}
\label{KoranyiVagitheorem} (Kor\'anyi-V\'agi type theorem) There is a constant $C>0$ such
that
\begin{equation*}
\int\limits_{\partial E} \vert M_{\rad}\mathcal{C}[f](\zeta)\vert^2 \vert d\zeta\vert \leq
C\int\limits_{\partial E} \vert f(\zeta)\vert^2 \vert d\zeta\vert
\end{equation*}
for every $f\in L^2(\partial E,\vert d\zeta\vert).$
\end{thm}

We recall the definition of the  Smirnov class $E^p,$ $p>0,$  on rectifiable
Jordan domains.
\begin{defi}
\label{defi6.8}
Let $p>0$ and $\Omega$  a  rectifiable Jordan domain. A function
$f\in\mathcal{O}(\Omega)$ is said to belong to the Smirnov class $E^p(\Omega)$ if
there exists a sequence of rectifiable closed Jordan curves
$(\mathcal{C}_n)_{n=1}^{\infty}$ in $\Omega,$ tending to the boundary
in the sense that  $\mathcal{C}_n$ eventually surrounds each compact
subdomain of $\Omega,$ such that
\begin{equation*}
\int\limits_{\mathcal{C}_n} \vert f(z)\vert^p \vert d z\vert \leq
 M<\infty,\qquad n\geq 1.
\end{equation*}
\end{defi}

Next, we rephrase some facts concerning the  Smirnov class $E^p,$ $p>0$  on rectifiable
Jordan domains in the context
  of  angular Jordan domains $\Omega(F,h).$
\begin{thm}
\label{Smirnovtheorem}
1) Let $\Omega$ be a rectifiable Jordan domain.
Then every $f\in E^p(\Omega)$ ($p>0$)  admits the  angular limit $f^{\ast}$ a.e.   on
$\partial\Omega.$\\
2) Let $\Omega:=\Omega(F,h)$ be an angular Jordan domain and let $\Gamma:=\partial\Omega.$
For any $0<\epsilon<\frac{h}{4},$ let $\Gamma_{\epsilon}$ be the
rectifiable closed Jordan curve defined in Part 4) of Proposition
\ref{prop6.2}.
Then $f\in E^p(\Omega)$ if  $\sup\limits_{0<\epsilon<\frac{h}{4}}
\int\limits_{\Gamma_{\epsilon}} \vert f(z)\vert^p \vert dz\vert<\infty.$
In addition,
 for an $f\in E^p(\Omega)$, $p>0,$ it holds that
\begin{equation*}
\int\limits_{\Gamma} \vert f^{\ast}(z)\vert^p \vert dz\vert \leq
\sup\limits_{0<\epsilon<\frac{h}{4}}
\int\limits_{\Gamma_{\epsilon}} \vert f(z)\vert^p \vert dz\vert.
\end{equation*}
3)  Every $f\in E^1(E)$ has a Cauchy representation
$f:= \mathcal{C}[f^{\ast}].$ Conversely, if $g\in L^1(\partial E,\vert
dz\vert)$ and
\begin{equation*}
\int\limits_{\partial E}  z^ng(z) d z=0, \qquad n=0,1,2,\ldots,
\end{equation*}
then $f:= \mathcal{C}[g]\in E^1(E)$
and $g$ coincides with $f^{\ast}$  a.e. on $\partial E.$
 \end{thm}
 \begin{proof} For the proof of Parts 1) and 3),
see \cite[p. 438--441]{go}. Taking into account Parts 6) and 7) of Proposition
\ref{prop6.2}, Part 2) also follows from the results in \cite[p.
438--441]{go}. Hence, the proof is complete.
 \end{proof}

\subsection{Proof of Theorem \ref{thm6.5}.  }
  We fix  $w_0\in G$  and  $0<\delta_0<\delta$ with
 $\omega(w_0,B,G)<\delta_0$ and
an arbitrary connected component $U$ of
$D_{\delta}:=\left\lbrace  z\in D:\ \omega(z,A,D) <1-\delta  \right\rbrace.$
Applying Theorem \ref{thm6.3}, we may find  an angular Jordan domain
$\Omega:=\Omega(F,h)\subset U$  such that $F\subset A\cap
A^{\ast}\cap U^{D}.$ In the course of the proof, the letter $C$
will denote a positive   constant that is not necessarily the same at each
 step.

Applying Carleman  Theorem (see, for example, \cite[p.2]{ai}), we have
\begin{eqnarray*}
f(z,b)&=&\lim\limits_{N\to\infty} \frac{1}{2\pi i}\int\limits_{A} e^{-N(g_1(a)-g_1(z))} \frac{f(a,b)da}{a-z},\qquad
z\in D,\ b\in B,\\
f(a,b)&=&\lim\limits_{r\to 1-}f(ra,b),\qquad a\in\partial D,\
b\in B.
\end{eqnarray*}
Consequently,  $f|_{\partial
D\times B}$ is measurable. In addition, by  (iii) this function
is   bounded.
Therefore, for every $N\in\N$ we are able to define  the function
$K_{\infty,N}(\cdot,w_0):\ \partial D\longrightarrow\C,$
\begin{equation}\label{eq6.4}
K_{\infty,N}(a,w_0):=\frac{1}{2\pi i}\int\limits_{B} e^{N(g_2(w_0)-g_2(b))}
\frac{f(a,b)db}{b-w_0},\qquad  a\in \partial D.
\end{equation}
Since, in virtue of (ii)--(iii), $f(a,\cdot)\in \mathcal{O}(G)$
and $\vert f(a,\cdot)\vert_{G}<C$ for $a\in A,$   it follows from Carleman Theorem
 that
\begin{equation}\label{eq6.5}
\lim\limits_{N\to\infty}K_{\infty,N}(a,w_0)=f(a,w_0),\qquad  a\in
A,
\end{equation}
and the above convergence is uniform with respect to $a\in A.$

 On the other hand,  by  (\ref{eq6.4}) we see that
 $K_{\infty,N}(\cdot,w_0)$ is measurable and
bounded. In addition, for any $n=0,1,2,\ldots,$ taking (ii) into account, we have that
\begin{equation*}
\int\limits_{\partial D}  K_{\infty,N}(a,w_0)a^nda
=\frac{1}{2\pi i}\int\limits_{B} \Big( \int\limits_{\partial D}f(a,b)a^nda\Big)
  \frac{ e^{N(g_2(w_0)-g_2(b))}db}{b-w_0}=0,
  \end{equation*}
where the first equality follows from an application of Fubini's Theorem and
the second one  from an application of Part 3) of Theorem \ref{Smirnovtheorem} to
$f(\cdot,b),$ $b\in B.$
Consequently, in virtue of Part 3) of  Theorem \ref{Smirnovtheorem}, we can  extend
 $K_{\infty,N}(\cdot,w_0)$ to $\overline{D}$ by setting
\begin{equation}\label{eq6.6}
K_{\infty,N}(z,w_0):= \mathcal{C}[K_{\infty,N}(\cdot,w_0) ](z)=\frac{1}{2\pi i}
\int\limits_{\partial D} \frac{K_{\infty,N}(a,w_0)da}{a-z},\qquad  z\in D.
\end{equation}
Then   the following identity holds
\begin{equation}\label{eq6.7}
\lim\limits_{z\to a,\ z\in\mathcal{A}_{\alpha}(a)}
K_{\infty,N}(z,w_0) =K_{\infty,N}(a,w_0)  ,\qquad
0<\alpha<\frac{\pi}{2},
\end{equation}
for a.e. $a\in \partial D.$

Now  we come back to the angular Jordan domain $\Omega$.  We keep the notation introduced in
Proposition \ref{prop6.2}.  For any $0<\epsilon<\frac{h}{4}$ and
any  $z\in \Gamma_{\epsilon},$ applying the Cauchy integral formula, we
obtain
\begin{multline}\label{eq6.8}
  K_{\infty,N}(z,w_0)-K_{N}(z,w_0)  \\
  =\frac{1}{(2\pi i)^2}
\int\limits_{\partial D\setminus A}\int\limits_{B}
e^{N(g_1(z)-g_1(a)) + N(g_2(w_0)-g_2(b))}
\frac{f(a,b)da db}{(a-z)(b-w_0)} \\
 =e^{N(g_1(z)-(1-\delta))}   \int\limits_{\partial D}\frac{p_N(a)da
}{a-z}.
\end{multline}
Using   the choice of $U$ and the hypothesis on $\delta$ and $\delta_0,$
it can be checked that
 \begin{equation}\label{eq6.9}
  \vert e^{N(g_1(\cdot)-(1-\delta))}\vert_U \leq 1, \qquad \vert p_N\vert_{\partial D}\leq
Ce^{-N(\delta-\delta_0)} .
\end{equation}
Therefore, recalling   the projection $\tau:\  E\setminus\{0\}\longrightarrow\partial E$
 (see Part 5) of Proposition \ref{prop6.2}), we estimate
\begin{multline}\label{eq6.10}
\int\limits_{\Gamma_{\epsilon}} \vert
K_{\infty,N}(z,w_0)-K_{N}(z,w_0)\vert^2 \vert dz\vert
\leq C\int\limits_{\Gamma_{\epsilon} }\Big\vert M_{\text{rad}}\mathcal{C}[ p_N](\tau(z))\Big\vert^2
\vert dz\vert\\
\leq 10 C\int\limits_{\tau\left(F_{0\epsilon}\cup F_{1\epsilon}\cup F_{2\epsilon}\right) }
\Big\vert M_{\text{rad}}\mathcal{C}[ p_N](a)\Big\vert^2
\vert da\vert
 +10 C\int\limits_{\tau (F_{3\epsilon} ) }
\Big\vert M_{\text{rad}}\mathcal{C}[ p_N](a)\Big\vert^2
\vert da\vert\\
\leq 20 C\int\limits_{ \partial E }
\Big\vert M_{\text{rad}}\mathcal{C}[ p_N](a)\Big\vert^2
\vert da\vert
\leq C\int\limits_{ \partial E } \vert p_N(a)\vert^2\vert
da\vert\leq Ce^{-N(\delta-\delta_0)} .
\end{multline}
  Here the first estimate follows from   (\ref{eq6.8})--(\ref{eq6.9}) and
the definition of the radial maximal function, the second and the third one are consequences of
 Part 5) of
Proposition \ref{prop6.2}, the fourth estimate holds by an application of
Theorem \ref{KoranyiVagitheorem}, and the last one follows from
(\ref{eq6.9}).

On the other hand, for any $0<\epsilon<\frac{h}{4},$
\begin{multline}\label{eq6.11}
\int\limits_{\Gamma_{\epsilon}} \vert
K_{N+1}(z,w_0)-K_{N}(z,w_0)\vert^2 \vert dz\vert \\
\leq 2\int\limits_{\Gamma_{\epsilon}} \vert
A_{N}(z,w_0) \vert^2 \vert dz\vert+ 2\int\limits_{\Gamma_{\epsilon}} \vert
B_{N}(z,w_0) \vert^2 \vert dz\vert
\leq Ce^{-N(\delta-\delta_0)},
\end{multline}
where $A_N$ and $B_N$ are given by formula (6) in \cite{go2}
and  the latter estimate follows from the same argument as in the proof of  (\ref{eq6.8})--(\ref{eq6.10}).
 We recall from (\ref{eq1_Gonchar_Carleman}) that
\begin{equation*}
\lim\limits_{N\to\infty}K_{N}(z,w_0)=K(z,w_0), \qquad
z\in\Gamma_{\epsilon}.
\end{equation*}
This, combined with (\ref{eq6.10})--(\ref{eq6.11}), implies that
\begin{equation}\label{eq6.12}
\int\limits_{\Gamma_{\epsilon}} \vert
K_{\infty,N}(z,w_0)-K(z,w_0)\vert^2 \vert dz\vert
 \leq C\cdot e^{-N(\delta-\delta_0)},\qquad 0<\epsilon<\frac{h}{4}.
\end{equation}
Since we have already shown that $\vert
K_{\infty,N}(\cdot,w_0)\vert_{D}<\infty,$ in virtue of Part 2) of
Theorem \ref{Smirnovtheorem}, we deduce from (\ref{eq6.12}) that $
K(\cdot,w_0)|_{\Omega}\in E^2(\Omega).$ For every $a\in \partial D,$
let $K(a,w_0)$ denote the angular limit of $K(\cdot,w_0)|_{\Omega}$
at $a$ (if the limit exists). It follows from  (\ref{eq6.12}) and
 Part 2) of
Theorem \ref{Smirnovtheorem} that
\begin{multline*}
\lim\limits_{N\to\infty}\int\limits_{\Gamma} \vert
K_{\infty,N}(a,w_0)-K(a,w_0)\vert^2 \vert da\vert\\
\leq \sup\limits_{0<\epsilon<\frac{h}{4}}
\int\limits_{\Gamma_{\epsilon}} \vert
K_{\infty,N}(z,w_0)-K(z,w_0)\vert^2 \vert dz\vert
 \leq \lim\limits_{N\to\infty} C\cdot e^{-N(\delta-\delta_0)}
=0.
\end{multline*}
This, combined with (\ref{eq6.5}) and Part 8) of Proposition \ref{prop6.2},  implies finally  that
\begin{equation*}
K(a,w_0)=f(a,w_0), \qquad\text{for a.e.}\ a\in F.
\end{equation*}
Hence, Theorem \ref{thm6.5} has been proved.
\hfill $\square$
\section{Proof of Theorem A for the case where $D$ and $ G$ are    Jordan domains}

Using an exhaustion argument, a compactness argument and  conformal mappings, the
case where  $D$ and $G$ are Jordan domains can be reduced to the following
case:

\smallskip

{\it  We assume that   $D=G=E,$ and  $\vert f\vert_W< 1.$ \hfill $(*)$}

\smallskip

Using hypotheses (i)--(iii) and $(*),$ we may apply Theorem \ref{thm6.5} and obtain a function
$K[f]
 \in\mathcal{O}(\widehat{W}^{\text{o}} ).$ Consequently, we are able to define the desired extension
 function $\hat{f}$ as follows
 \begin{equation*}
 \hat{f}:=K[f].
 \end{equation*}
In this section we   will use repeatedly  Part 3) of Theorem \ref{thm5.15}
\begin{equation*}
\omega(\cdot, A,\Omega)=\omega(\cdot, A^{\ast},\Omega),
\end{equation*}
where $\Omega\subset \C$ is an open set and $A$ is a Jordan measurable subset of $\partial \Omega.$

  The remaining part of the proof is divided into
several steps.

\smallskip

\noindent{\bf Step 1:} {\it  Proof of the estimate
\begin{equation*}
\vert \hat{f}\vert_{\widehat{W}^{\text{o}}} \leq \vert f\vert_W.
\end{equation*}
}

\smallskip

 \noindent{\it Proof of Step 1.}
 Let $(z_0,w_0)$ be an  arbitrary point of $ \widehat{W}^{\text{o}}.$ Then we may find
an $\delta\in (0,1)$ such that $0<\omega(w_0,B,G) <\delta<1-\omega(z_0,A,D) .$
 Let $U$ be the  connected component  of
$D_{\delta}:=\left\lbrace  z\in D:\ \omega(z,A,D) <1-\delta  \right\rbrace$
 that contains $z_0$.
By Theorem \ref{thm6.3} we may find  an angular Jordan domain
$\Omega:=\Omega(F,h)$    contained in $U$  such that $F\subset A\cap
A^{\ast}\cap U^{D}.$ In addition,
 for every $N\in \N,$ applying Theorem \ref{thm6.5} to the function $ f^N,$
 we obtain the function $K[f^N]\in\mathcal{O}(\widehat{W}^{\text{o}})$
 with the following property
\begin{multline*}
\lim\limits_{z\to a,\ z\in\mathcal{A}_{\alpha}(a)}K[f^N](z,w_0 )=f(a,w_0
)^N\\
=\lim\limits_{z\to a,\ z\in\mathcal{A}_{\alpha}(a)}\left(K[f](z,w_0 )\right)^N,\
 0< \alpha<\frac{\pi}{2},
\end{multline*}
for a.e. $a\in F.$

 Consequently, an application of  Theorem \ref{thm6.4}  gives that
 \begin{equation*}
K[f^N](z_0,w_0)=\left(K[f](z_0,w_0)\right)^N,\ N\in
\N,
\end{equation*}
 Since $(z_0,w_0)\in\widehat{W}^{\text{o}}$ is arbitrarily chosen,   it follows from the latter identity that
\begin{equation}\label{eq7.2}
K[f^N](z,w)=\left(K[f](z,w)\right)^N,\ N\in \N,\ (z,w)\in\widehat{W}^{\text{o}}.
\end{equation}

 Now we are able to  conclude the
proof in the same way as in \cite[p. 23]{go2}. More precisely,
taking into account  (\ref{eq7.2}), one
gets that
\begin{equation*}
\vert \hat{f}^N(z,w)\vert \leq\vert K[f^N](z,w)\vert \leq \frac{C\vert f\vert_W^N}{(1-\vert z\vert)  (1-\vert  w\vert)
   (1-e^{-(1-\omega(z,w))})},\
 (z,w)\in\widehat{W}^{\text{o}}.
\end{equation*}
Extracting the $N$th roots of both sides and letting $N$ tend to $\infty,$
the desired estimate of Step 1 follows. \hfill $\square$

\smallskip

\noindent{\bf Step 2:} {\it We shall prove that $\hat{f}$ is the unique function $\mathcal{O}(
\widehat{W}^{\text{o}})$
 which verifies  Property 1).}

\smallskip

 \noindent{\it Proof of Step 2.} First we show that the
 function $\hat{f} $  satisfies Property 1).
 Without loss of generality, it suffices to prove that there is a subset $\tilde{B}$ of
 $B\cap B^{\ast}$
such that $\mes(\tilde{B})=\mes(B)$ and $\hat{f}$ admits the  angular limit $f$ at
every point of $D\times \tilde{B} .$

For any $a\in A$ put
\begin{equation*}
B_{a}:=\left\lbrace b\in B:\ f(a,\cdot)\ \text{has an angular limit at}\
b\right\rbrace.
\end{equation*}
By hypothesis (iii), we have
$\mes(B_{a})=\mes(B),$ $ a\in A.$  
Consequently, applying Fubini's Theorem, we obtain that
\begin{equation*}
\int\limits_{A}\mes(B_{a})\vert d a\vert=\mes(A)\mes(B)
=\int\limits_{B}\mes\left(\left\lbrace a\in A:\ b\in B_{a}\right\rbrace\right)\vert
d b\vert.
\end{equation*}
Hence,
\begin{equation}\label{eq7.2.1}
 \mes\left(\left\lbrace a\in A:\ b\in
B_{a}\right\rbrace\right)=\mes(A)\qquad\text{ for a.e.}\ b\in B.
\end{equation}
The same reasoning also gives that
\begin{equation}\label{eq7.2.2}
\mes\left(\left\lbrace a\in A:\ f(a,b)=
 f_1(a,b)\right\rbrace\right)=\mes(A)\qquad\text{ for a.e.}\ b\in B.
\end{equation}
Set
\begin{multline}\label{eq7.2.3}
\tilde{B}:=\left\lbrace b\in B\cap B^{\ast}:\
\mes\left(\left\lbrace a\in A:\ b\in
B_{a}\right\rbrace\right)=\mes(A)\right.\\
\left. \text{ and}\  \mes\left(\left\lbrace a\in A:\
f(a,b)= f_1(a,b)\right\rbrace\right)=\mes(A)
\right\rbrace. \end{multline}
We deduce from (\ref{eq7.2.1})--(\ref{eq7.2.3}) that
\begin{equation}\label{eq7.2.4}
 \mes(\tilde{B})=\mes(B).
 \end{equation}

Fix an arbitrary point $b_0\in \tilde{B}$ and
let $(w_n)_{n=1}^{\infty}$ be
an arbitrary sequence of $G$ such that $\lim\limits_{n\to\infty} w_n=b_0$
and $w_n\in\mathcal{A}_{\alpha}(b_0)$ for some fixed number  $0<\alpha<\frac{\pi}{2}.$
Fix an arbitrary point $z_0$ of $D$ and let  $(z_n)_{n=1}^{\infty}$ be
an arbitrary sequence of $D$ such that $\lim\limits_{n\to\infty}z_n=z_0.$

 Clearly, we may find  $0<\delta_1<1$ such that
\begin{equation}\label{eq7.2.5}
 \sup\limits_{n\in\N}\omega(z_n,A,D) <1-\delta_1.
\end{equation}
Fix  an $\delta_2$ such that $0<\delta_2<\delta_1 .$
Since $b_0$ is locally regular relative to $B$ and
$\lim\limits_{n\to\infty} w_n=b_0$ and
$w_n\in\mathcal{A}_{\alpha}(b_0),$
there is a sufficiently large number $N_0$ with
\begin{equation}\label{eq7.2.6}
 \omega(w_{n},B,G)<\delta_2,\qquad n> N_0.
\end{equation}
Let $U$ be that connected component  of the
 following open set
 \begin{equation*}
D_{\delta_1}:=\left\lbrace  z\in D:\ \omega(z,A,D) <1-\delta_1 \right\rbrace
\end{equation*}
 which contains $z_0$  (see (\ref{eq7.2.5})).
Applying Theorem \ref{thm6.3}, we may find  an angular Jordan domain
$\Omega:=\Omega(F,h)$    contained in $U$  such that $F\subset A\cap
A^{\ast}\cap U^{D}.$
Let $V$  be a rectifiable Jordan domain with   $\Omega\subset V\subset  U,$
 $w_0\in V$, and
 $V\cap \mathcal{U}=\Omega\cap \mathcal{U}$ for some
   neighborhood  $\mathcal{U}$  of the base $F$ of $\Omega.$

 In virtue of (\ref{eq7.2.6}) and of the fact that $V\subset U\subset D_{\delta}
 ,$ we obtain that
\begin{equation}\label{eq7.2.7}
V\times \{w_n\} \subset
\widehat{W}^{\text{o}},\qquad n> N_0.
\end{equation}
Consequently,   Theorem \ref{thm6.5} yields that for any $ n>N_0,$
\begin{equation}\label{eq7.2.8}
f(a,w_n)=\lim\limits_{z\to a,\ z\in\mathcal{A}_{\alpha}(a)}
\hat{f}(z,w_n),\qquad   0< \alpha<\frac{\pi}{2},
\end{equation}
for a.e. $a\in F.$

Next, for any $n> N_0$ let
\begin{equation*}
\begin{split}
F_n&:=\left\lbrace a\in F:\  b_0\in B_{a}\ \text{and}\ f(a,w_n)
=\lim\limits_{z\to a,\
z\in\mathcal{A}_{\alpha}(a)}
\hat{f}(z,w_n)\right\rbrace,\\
F_0&:=\bigcap\limits_{n=N_0+1}^{\infty} F_n.
\end{split}
\end{equation*}
It follows from (\ref{eq7.2.3}), (\ref{eq7.2.8}) and the fact that $b_0\in \tilde{B}$
that $\mes(F_n)=\mes(F),$
$n>N_0.$ Hence
\begin{equation}\label{eq7.2.9}
\mes(F_0)=\mes(F)>0.
\end{equation}

In virtue of  (\ref{eq7.2.7}), consider the following   holomorphic functions on $V$
\begin{equation}\label{eq7.2.10}
h_n(t):=\hat{f}(t,w_n) \quad \text{and}\quad
h_0(t):=f(t,b_0),\qquad  t\in V,\ n> N_0 .
\end{equation}
 Since we have already shown in Step I that  $\vert h_n\vert_{V}\leq \vert f\vert_X <\infty,$ $n>
 N_0$ or $n=0,$
  applying Part 1) of Theorem \ref{Smirnovtheorem},
 we may find a subset $\Delta$ of $F_0$ with
$\mes(\Delta)=\mes(F_0)>0$ such that   $h_n,$ $n>N_0$ (resp. $h_0$) admits
the angular limit
$f_1\left(t,w_n\right)$ (resp. $f_1\left(t,b_0\right))$ at $t\in \Delta.$
 Observe that by (\ref{eq7.2.3}) and the fact that   $b_0\in
 \tilde{B}$ we have that
\begin{equation*}
\lim\limits_{n\to\infty} f_1\left(t,w_n\right)=f_1\left(t,b_0\right)=f\left(t,b_0\right)
\qquad \text{for a.e.}\ t\in \Delta.
\end{equation*}

   Using this and  (\ref{eq7.2.10}), we are able to apply  Khinchin--Ostrowski Theorem  (see  \cite[Theorem 4, p. 397]{go})
   to the sequence $(h_n)_{n=0}^{\infty}.$ Consequently, one gets
\begin{equation*}
 \lim\limits_{n\to\infty} \hat{f}(z_n,w_n)=f(z_0,b_0).
\end{equation*}
This shows that $\hat{f}$ admits the  angular limit $f$ at
every point of $ D\times\tilde{B} .$ Hence, $\hat{f}$ satisfies Property 1).

In order  to complete Step 2  we need to show the uniqueness of $\hat{f}$. To do this, let
$\hat{\hat{f}}\in  \mathcal{O}(\widehat{W}^{\text{o}})$ be  a function with the following property:
There is a subset $\tilde{\tilde{A}}$   (resp.   $\tilde{\tilde{B}}$) of
 $A\cap A^{\ast}$   (resp.   $b\cap B^{\ast}$)
such that $\mes(A\setminus\tilde{ \tilde{A}})=\mes(B\setminus\tilde{ \tilde{B}})  =       0$
 and $\hat{\hat{f}}$ admits the  angular limit $f$ at
every point of $(\tilde{\tilde{A}}\times G)\cup (D\times\tilde{\tilde{B}}) .$
Fix an arbitrary point $(z_0,w_0)\in \widehat{W}^{\text{o}}.$  Let $U$ be the connected component
containing $z_0$ of the following open set
\begin{equation*}
\left\lbrace z\in D:\ \omega(z,A,D)<1-\omega(w_0,B,G)
\right\rbrace.
\end{equation*}
We deduce from the  property of $\hat{f}$ and $\hat{\hat{f}}$ that both
holomorphic functions $\hat{f}(\cdot,w_0)|_U$ and $\hat{\hat{f}}(\cdot,w_0)|_U$
 admit the angular limit $f(\cdot,w_0)$ at every point of
 $\tilde{A}\cap\tilde{\tilde{A}}\cap U^{D}.$ Consequently, applying
 Theorem \ref{thm6.4} yields that
 $\hat{f}(\cdot,w_0)=\hat{\hat{f}}(\cdot,w_0)$ on $U.$ Hence,
  $\hat{f}( z_0,w_0)=\hat{\hat{f}}(z_0,w_0).$ Since $(z_0,w_0)\in\widehat{W}^{\text{o}}$
  is arbitrary, the uniqueness of $\hat{f}$ is established.
This completes Step 2.
\hfill $\square$

\smallskip

\noindent{\bf Step 3:} {\it   Proof of Part 2).
}

\smallskip

 \noindent{\it Proof of Step 3.}
  Fix  $(z_0,w_0)\in\widehat{W}^{\text{o}}.$
For every  $b\in B$ we have
\begin{equation*}
\vert f(a,b)\vert \leq \vert f\vert_{A\times
B},\ a\in A, \qquad\text{and}\qquad \vert f(z,b)\vert \leq \vert f\vert_{W},\ z\in
D.
\end{equation*}
Therefore, the Two-Constant Theorem (see Theorem 2.2 in \cite{pn}) implies that
\begin{equation}\label{eq7.3.1}
\vert f(z,b)\vert \leq \vert f\vert_{A\times
B}^{1-\omega(z,A,D)} \vert f\vert_{W}^{\omega(z,A,D)},\qquad z\in D,\
b\in B.
\end{equation}
 Let $\delta:=\omega(z_0,A,D)$ and consider the    $\delta$-level set
 \begin{equation*}
G_{\delta}:=\left\lbrace  w\in G: \omega(w,B,G)<1-
\delta\right\rbrace.
 \end{equation*}
 Clearly,  $w_0\in G_{\delta}.$

 Recall from  Step 2   that
$\tilde{B}\subset B\cap B^{\ast},$  $\mes\Big( (B\cap B^{\ast})\setminus
\tilde{B}\Big)=0, $   and
\begin{equation}\label{eq7.3.2}
f(z_0,b)=\lim\limits_{w\to b,\
w\in\mathcal{A}_{\alpha}(b)}
\hat{f}(z_0,w),\qquad  0<\alpha<\frac{\pi}{2},\ b\in \tilde{B}.
\end{equation}

 Consider the following  function $h:\  G_{\delta}\cup
 \tilde{B}\longrightarrow\C$ defined  by
 \begin{equation}\label{eq7.3.3}
 h(t):=
\begin{cases}
 \hat{f}(z_0,t),
  & t\in   G_{\delta} \\
 f(z_0,t), & t \in \tilde{B}
\end{cases}.
\end{equation}
Clearly, $h|_{ G_{\delta}}\in\mathcal{O}(G_{\delta}  ).$

On the other hand, in virtue of  (\ref{eq7.3.3}) and  the result of Step 1, we
have
 \begin{equation} \label{eq7.3.4}
 \vert h\vert_{G_{\delta}}
 \leq \vert \hat{f} \vert_{\widehat{W}^{\text{0}}}\leq \vert f\vert_{W}<\infty.
  \end{equation}
 In addition, applying
   Corollary \ref{cor5.21}  and taking (\ref{eq7.3.2})--(\ref{eq7.3.3}) into account yields
 \begin{equation*}
  \vert h(t)\vert
  \leq \vert h\vert_{ \tilde{B}}^{1-\omega_{\delta}(t,A,D)} \vert
 h\vert_{ G_{\delta}}^{\omega_{\delta}(t,A,D)},\quad   t\in G_{\delta} ,
  \end{equation*}
 where,  by Theorem \ref{thm5.20},
\begin{equation*}
\omega_{\delta}(t,B ,G)=\frac{\omega(t,B,G)}{1-\omega(z_0,A, D)}.
\end{equation*}
 This,   combined with  (\ref{eq7.3.1})--(\ref{eq7.3.4}), implies that
\begin{eqnarray*}
\vert \hat{f}(z_0,w_0)\vert &=&\vert h(w_0)\vert\leq
 \vert f\vert_{A\times B}^{1-\omega(z_0,A,D)-\omega(w_0,B,G)} \vert
 f\vert_{W}^{\omega(z_0,A,D)+\omega(w_0,B,G)}.
 \end{eqnarray*}
Hence Part 2) for the point $(z_0,w_0)$ is proved.
\hfill  $\square$

\smallskip

\noindent{\bf Step 4:} {\it   Proof of Part 3).
}

\smallskip

 \noindent{\it Proof of Step 4.}
Let $(a_0,w_0)\in A^{\ast}\times G$ be
such that the following limit exists
\begin{equation*}
\lambda:= \lim\limits_{(a,w)\to (a_0,w_0),\ (a,w)\in A\times G}
f(a,w).
 \end{equation*}
 We like to show that
$\hat{f}$ admits the angular limit $\lambda$ at  $(a_0,w_0).$

  For any $0<\epsilon<\frac{1}{2},$
 we may find an open neighborhood  $A_{a_0}$ of $a_0$ in $A$
and a positive number $r>0$
such that   $\B(w_0,r)\Subset G$  and
\begin{equation}\label{eq7.4.1}
\left\vert f(a,w)-\lambda\right\vert <\epsilon^2,\qquad
a\in A_{a_0},\   \vert w-w_0\vert \leq r.
\end{equation}
Put
\begin{equation}\label{eq7.4.2}
\delta:=\sup\limits_{w\in \overline{\B(w_0,r)}} \omega(w,B,G).
\end{equation}
Since $a_0\in A^{\ast},$   it is clear  that  $\mes(A_{a_0} )>0.$ Next, consider
the level set
\begin{equation*}
D_{\delta}:=\left\lbrace z\in D:\ \omega(z,A_{a_0},D)<1-\delta
\right\rbrace.
\end{equation*}
In virtue of (\ref{eq7.4.2}), we can define
\begin{equation}\label{eq7.4.3}
h(t,w):=\hat{f}(t,w)-\lambda ,\qquad  t\in D_{\delta},\ w\in \overline{\B(w_0,r)}.
\end{equation}
Clearly,
 \begin{equation}\label{eq7.4.4}
  \vert h\vert_{D_{\delta}}
 \leq 2\vert \hat{f} \vert_{\widehat{W}^{\text{0}}}= 2 \vert f\vert_{W}=2.
  \end{equation}
By   (\ref{eq7.4.3}) and using the result of  Step 2,   we know
that for every $w\in \B(w_0,r)$
the holomorphic function $h(\cdot,w)|_{D_{\delta}}$  admits the angular limit
$f(a,w)-\lambda$ at $a$ for  $a\in\tilde{A}\cap A_{a_0},$
where  $\tilde{A}$ is given in Step 2.
Consequently, applying  Corollary \ref{cor5.21} and taking (\ref{eq7.4.1}) and (\ref{eq7.4.4})
into account, we see that
\begin{equation*}
\left\vert h(t,w) \right\vert
<\epsilon^{2(1-\omega_{\delta}(t,A_{a_0},D))}2^{\omega_{\delta}(t,A_{a_0},D)},\qquad
t\in D_{\delta}.
\end{equation*}
Let  $0<\alpha<\frac{\pi}{2}.$
 In virtue of Theorem \ref{thm5.20} and the hypothesis that $a_0\in A^{\ast},$ we  deduce  that
 $\lim\limits_{t\to a_0,\ t\in \mathcal{A}_{\alpha}(a_0)}
 \omega_{\delta}(t,A_{a_0},D)=0.$ Consequently,  there is
 an $r_{\alpha}>0$ such that
\begin{equation*}
\left\vert f(z,w)-\lambda\right\vert= \left\vert h(z,w) \right\vert<\epsilon,\qquad
z\in\mathcal{A}_{\alpha}(a_0)\cap \{\vert z-a_0\vert <r_{\alpha} \},\ w\in
\B(w_0,r).
\end{equation*}
This completes the above assertion.

 Similarly, we can prove that
$\hat{f}$ admits the angular limit
\begin{equation*}
\lim\limits_{(z,b)\to (z_0,b_0),\ (z,b)\in D\times B}
f(z,b)
\end{equation*}
 at any point $(z_0,b_0),$ if the latter limit exists.
Hence the proof of  Step 4 (i.e. Part 3)) is finished.  \hfill  $\square$

\smallskip

\noindent{\bf Step 5:} {\it   Proof of Part 4).
}

\smallskip

 \noindent{\it Proof of Step 5.}
 Let $(a_0,b_0)\in A^{\ast}\times B^{\ast}$ be
such that the following limit exists
\begin{equation*}
\lambda:= \lim\limits_{(a,b)\to (a_0,b_0),\ (a,b)\in A\times B}
 f(a,b).
 \end{equation*}
 We like to show that
$\hat{f}$ admits the angular limit $ \lambda$ at  $(a_0,b_0).$

 Recall that $\vert f\vert_X<1,$   and fix  an arbitrary $0<\epsilon<\frac{1}{2}.$
  Since $(a_0,b_0)\in A^{\ast}\times B^{\ast},$ we may
 find  an open neighborhood $A_{a_0}$ of $a_0$ in $A$ (resp.
an open neighborhood $B_{b_0}$ of $b_0$ in $B)$  such that
\begin{equation}\label{eq7.5.1}
\left\vert f(a,b)-\lambda\right\vert <\epsilon^2,\qquad
a\in A_{a_0}   ,\ b\in B_{b_0}.
\end{equation}
 It is clear that
  $\mes(A_{a_0})>0$ and  $\mes(B_{b_0} )>0.$

  Consider the function
\begin{equation}\label{eq7.5.2}
h(z,w):=f(z,w)-\lambda ,\qquad  (z,w)\in \X(A_{a_0},B_{b_0}; D,G).
\end{equation}
Clearly,
\begin{equation}\label{eq7.5.3}
\vert h(z,w) \vert \leq 2,\qquad  (z,w)\in \X(A_{a_0},B_{b_0}; D,G).
\end{equation}
Applying  the results of Steps 1--3  to $h,$ we obtain the function
\begin{equation}\label{eq7.5.4}
\hat{h}:=K[h]\qquad\text{on}\ \ \widehat{\X}^{\text{0}}(A_{a_0},B_{b_0}; D,G).
\end{equation}
so that $\hat{h}$ admits the angular limit $h$ on $(\tilde{A}_{a_0}\times G)
\cup (D\times \tilde{B}_{b_0}),$ where  $\tilde{A}_{a_0}, $ $\tilde{B}_{b_0}$ are
given by Step 2. Clearly,
\begin{equation*}
  \widehat{\X}(A_{a_0},B_{b_0}; D,G)\subset  \widehat{\X}(A,B; D,G).
\end{equation*}
Consequently, arguing as in Step 1 and taking into account the above mentioned angular limit of $\hat{h},$
 we conclude that
\begin{equation*}
\hat{h}=\hat{f}-\lambda\qquad\text{on}\ \widehat{\X}(A_{a_0},B_{b_0}; D,G).
\end{equation*}
Consequently, applying Step 3 and taking into account (\ref{eq7.5.1})--(\ref{eq7.5.4})
and the inequality $\vert f\vert_X<1,$  we  see that
\begin{equation*}\begin{split}
\left\vert \hat{f}(z,w)-\lambda\right\vert&= \vert \hat{h}(z,w)\vert \leq
\vert h\vert_{A_{a_0}\times B_{b_0}}^{1- \omega(z,A_{a_0},D)- \omega(w,B_{b_0},G)}
(2\vert f\vert_X)^{ \omega(z,A_{a_0},D)+ \omega(w,B_{b_0},G) }\\
&< \epsilon^{2\Big(1-\omega(z,A_{a_0},D)- \omega(w,B_{b_0},G)   \Big )}2^{
 \omega(z,A_{a_0},D)+ \omega(w,B_{b_0},G)}.
\end{split}
\end{equation*}
Therefore, for  all $(z,w)\in\widehat{\X}(A_{a_0},B_{b_0}; D,G)$ satisfying
\begin{equation}\label{eq7.5.5}
\omega(z,A_{a_0},D)+ \omega(w,B_{b_0},G)<\frac{1}{3},
\end{equation}
we deduce from the latter estimate that
\begin{equation}\label{eq7.5.6}
\left\vert \hat{f}(z,w)-\lambda\right\vert<\epsilon.
\end{equation}
 Since $a_0$ (resp. $b_0$  ) is locally regular relative to $A_{a_0}$ (resp.
 $B_{b_0}$),
 there is an $r_{\alpha}>0$ such that  (\ref{eq7.5.5}) is fulfilled for
\begin{equation*}
(z,w)\in \left(\mathcal{A}_{\alpha}(a_0)\cap \{\vert z-a_0\vert <r_{\alpha}
\}\right) \times\left(
 \mathcal{A}_{\alpha}(b_0)\cap \{\vert w-b_0\vert <r_{\alpha} \}\right).
\end{equation*}
This, combined with (\ref{eq7.5.6}), completes the proof.   Hence Step 5 (i.e. Part 4)) is finished.
 \hfill $\square$

\smallskip

\noindent{\bf Step 6:} {\it   Proof of Part 5).
}

\smallskip

 \noindent{\it Proof of Step 6.}
 In virtue of Step 5, we only need to show that $\hat{f}$ admits the angular
 limit  $f$ on $(A^{\ast}\times G)\cup (D\times B^{\ast}).$
  To do this let $(a_0,w_0)\in A^{\ast}\times G$
 and choose an arbitrary $0<\epsilon< 1.$ Fix  a    compact subset $K$ of $B\cap B^{\ast}$
 such that $\mes(K)>0$ and a sufficiently large $N$ such that
 \begin{equation} \label{eq7.6.1}
\epsilon^{N(1-\omega(w_0,K,G))}(2\vert f\vert_X)^{\omega(w_0,K,G)}<\frac{\epsilon}{2}.
 \end{equation}
Using the hypothesis that $f$ can be extended to a continuous function on $A^{\ast}\times
B^{\ast},$ we may find  an open neighborhood $ A_{a_0}$
 of $a_0$ in $A^{\ast}$ such that
\begin{equation} \label{eq7.6.2}
\left\vert f(a,b)-f(a_0,b) \right\vert\leq \epsilon^N,\qquad
a\in A_{a_0}\cap A^{\ast}_{a_0},\ b\in K .
\end{equation}
On the other hand,
\begin{equation}\label{eq7.6.3}
\left\vert f(a,w)-f(a_0,w_0) \right\vert\leq 2\vert f\vert_X< 2,\qquad
  a\in A_{a_0}\cap A^{\ast}_{a_0},\ w\in G.
\end{equation}
For  $a\in A_{a_0}\cap A^{\ast}_{a_0}, $ applying the Two-Constant Theorem  to the function
$f(a,\cdot)-f(a_0,\cdot)\in\mathcal{O}(G)$ and taking (\ref{eq7.6.1})--(\ref{eq7.6.3}) into
account,   we deduce that
\begin{equation}\label{eq7.6.4}
\left\vert f(a,w_0)-f(a_0,w_0) \right\vert\leq \epsilon^{N(1-\omega( w_0,K,G  ))}
(2\vert f\vert_X)^{\omega(w_0,K,G) }<\frac{\epsilon}{2}.
\end{equation}
  Since $  f(a,\cdot)|_{G} $ is a bounded holomorphic
function for $a\in A,$ there is an open  neighborhood $V$ of $w_0$ such that
\begin{equation*}
\left\vert f(a,w)-f(a,w_0) \right\vert  <\frac{\epsilon}{2},\qquad
a\in A,\ w\in V.
\end{equation*}
This, combined with (\ref{eq7.6.4}), implies that
\begin{eqnarray*}
\left\vert f(a,w)-f(a_0,w_0) \right\vert &\leq& \left\vert f(a,w_0)-f(a_0,w_0) \right\vert
+\left\vert f(a,w)-f(a,w_0) \right\vert  \\
&<&\frac{\epsilon}{2}+ \frac{\epsilon}{2}=\epsilon ,\qquad
a\in A_{a_0},\ w\in V.
\end{eqnarray*}
 Therefore, $f$ is continuous at  $(a_0,w_0).$ Consequently, we conclude, by
 Step 4, that
$\hat{f}$ admits the angular limit $f(a_0,w_0)$ at  $(a_0,w_0).$
Similarly, we may also show that $\hat{f}$ admits the angular limit $f(z_0,b_0)$
at   every point  $(z_0,b_0)\in D\times B^{\ast}.$
This completes the proof of  the last step. \hfill $\square$
\section{Preparatory results}

We first develop some auxiliary results.  This preparation will enable us to generalize the results
of section 6 to the general case considered in  Theorem A.

\begin{defi}\label{defi8.1}
Let $\Omega$ be a  complex manifold of dimension $ 1$  and  $A\subset \Omega.$ Define
\begin{equation*}
\omega(\cdot,A,\Omega):=\sup\left\lbrace u:\ u\in\mathcal{SH}(\Omega),\ u\leq 1\ \text{on}\
\Omega, \ u\leq 0\ \text{on}\ A
  \right\rbrace.
\end{equation*}
The function $\omega(\cdot,A,\Omega)$ is called the {\it
the harmonic measure of $A$ relative to $\Omega.$}
A point $\zeta\in \overline{A}\cap \Omega$ is said  to be a {\it locally regular point relative to $A$}
 if  \begin{equation*}
\lim\limits_{z\to \zeta}  \omega(z,A\cap U,\Omega\cap U)=0
\end{equation*}
for any open neighborhood $U$ of $\zeta.$
If, moreover, $\zeta\in A,$ then $\zeta$ is said to be
a   {\it locally regular point of  $A.$}
The set of all locally regular points relative to $A$ is denoted by $A^{\ast}.$
$A$ is said to be locally regular if  $A=A^{\ast}.$
\end{defi}

\begin{prop}\label{prop8.2}
 Let $X$ be a complex manifold of dimension $1,$
 $D\subset X$  an open set and  $A\subset \partial D$   a Jordan  measurable
subset  of positive length.
Let $\{a_j\}_{j\in J}$ be a finite or countable
subset of $A$ with the following properties:
\begin{itemize}
\item[(i)]
 For any $j\in J,$  there is an
open neighborhood $U_{j}$ of  $a_j$ such that $D\cap U_{j}$ is   either a
 Jordan domain or the disjoint union of two Jordan
domains;
\item[(ii)] $A\subset \bigcup\limits_{j\in J} U_j.$
\end{itemize}
For any $0<\delta<\frac{1}{2},$   define
\begin{equation*}
\begin{split}
U_{j,\delta}&:=\left\lbrace z\in  D\cap U_{j}:\ \omega(z, A^{\ast}\cap U_j, D\cap U_j)<\delta  \right\rbrace,\qquad
j\in J,\\
A_{\delta}&:=\bigcup\limits_{j\in J} U_{j,\delta},\\
D_{\delta}&:=\left\lbrace z\in D:\ \omega(z, A^{\ast},  D)<1-\delta  \right\rbrace.
\end{split}
\end{equation*}
Then:\\
1) $ A\cap A^{\ast}\subset A_{\delta}^D$ and $A_{\delta}\subset D_{1-\delta} \subset D_{\delta};$\\
2) $\omega(z,A^{\ast},D)-\delta
\leq \omega(z,A_{\delta},D) \leq\omega(z,A^{\ast},D),$  $ z\in D .$
\end{prop}
\begin{proof} To prove Part 1), let $a\in A\cap A^{\ast}$ and fix an $j\in J$ such that $a\in U_j.$ Then
\begin{equation*}
\lim\limits_{z\to a,\ z\in\mathcal{A}_{\alpha}(a)} \omega(z, A^{\ast}\cap U_j, D\cap
U_j)=0,\qquad 0<\alpha<\frac{\pi}{2}.
\end{equation*}
Consequently, for every $ 0<\alpha<\frac{\pi}{2},$ there is an open
neighborhood $V_{\alpha}\subset U_j$ of $a$ such that
\begin{equation*}
  \omega(z, A^{\ast}\cap U_j, D\cap
U_j)<\delta,\qquad  z\in \mathcal{A}_{\alpha}(a)\cap V_{\alpha}.
\end{equation*}
 This proves $ A\cap A^{\ast}\subset A_{\delta}^D.$

To prove the second assertion of Part 1), one applies the Subordination Principle
and obtains for $z\in U_{j,\delta},$
\begin{equation}\label{eq8.1.1}
\omega(z,A^{\ast},D)\leq \omega(z,A^{\ast}\cap U_j,D\cap U_j)
 <\delta<1-\delta.
\end{equation}
Hence, $z\in D_{1-\delta}.$  This implies that   $A_{\delta}\subset
D_{1-\delta}.$ In addition, since  $0<\delta<\frac{1}{2},$ it follows that
 $ D_{1-\delta}\subset D_{\delta}.$ Hence,  Part 1) is proved.

We turn to Part 2). Since $A_{\delta}$ is an open set and, by Part 1),
$ A\cap A^{\ast}\subset A_{\delta}^D,$  it follows from Definitions
\ref{endpoint}, \ref{defi8.1}   that
 \begin{equation*}
\omega(z,A_{\delta},D)
\leq  \omega(z,A\cap A^{\ast},D),\qquad z\in D.
\end{equation*}
Hence, in virtue of  Theorem  \ref{thm5.15}, it follows
that
\begin{equation*}
 \omega(z,A_{\delta},D) \leq  \omega(z,A^{\ast} ,D),\qquad z\in D,
\end{equation*}
which proves the second  estimate of Part 2).

To complete Part 2), let $z\in A_{\delta}.$  Choose  $j\in J$
 such that $z\in U_{j,\delta}.$
 We deduce from (\ref{eq8.1.1})  that
$\omega(z,A^{\ast} ,D)-\delta\leq 0.$   Hence,
\begin{equation*}
\omega(z,A^{\ast} ,D)-\delta\leq 0, \qquad z\in
A_{\delta}.
\end{equation*}
On the other hand,  $\omega(z,A^{\ast} ,D)-\delta<1,$ $z\in D.$
Consequently, the first estimate  of Part 2)
follows. The proof of the  lemma is   finished.
\end{proof}

\smallskip

The main ingredient in the proof of Theorem A  is the following  mixed cross theorem.
\begin{thm}\label{mixedcrossthm}
Let $X$  and $Y$ be complex manifolds of dimension $1,$
 $D\subset X$ and  $\Omega\subset Y$   open subsets, and
  $A\subset  D$ and $  B\subset \partial \Omega.$
  Assume that
  $A=\bigcup\limits_{k=1}^{\infty} A_k$
  with $A_k$ locally regular  
   compact subsets of $D,$
  $A_k\subset A_{k+1},$ $k\geq 1.$  In addition,
  $B\subset \partial \Omega$  is a Jordan measurable subset of positive length. For $0\leq\delta<1$ put
    $G:=\left\lbrace
   w\in \Omega:\ \omega(w,B,\Omega)<1-\delta\right\rbrace.$
    Let $W:= \X(A,B;D,G)$, $W^{\text{o}}:= \X^{\text{o}}(A,B;D,G),$
      and   (using the notation $\omega_{\delta}(\cdot,B,\Omega)$ of Theorem \ref{thm5.20})   \begin{equation*}
    \widehat{W}^{\text{o}}=
     \widehat{\X}^{\text{o}}(A,B;D,G):=\left\lbrace (z,w)\in D\times G:\  \omega(z,A^{\ast},D)+\omega_{\delta}
     (w,B,\Omega)<1 \right\rbrace.
     \end{equation*}
Let  $f:\ W\longrightarrow \C$ be   such that
\begin{itemize}
\item[(i)] $f\in\mathcal{O}_s(W^{\text{o}});$
\item[(ii)]  $f$ is  Jordan measurable and locally  bounded on $W;$
\item[(iii)] for any $z\in A,$
\begin{equation*}
\lim\limits_{w\to \eta,\
w\in\mathcal{A}_{\alpha}(\eta) }f(z,w)=f(z,\eta),\qquad \eta\in B,\
0<\alpha<\frac{\pi}{2}.
\end{equation*}
\end{itemize}
      Then there is a unique function
$\hat{f}\in\mathcal{O}(\widehat{W}^{\text{o}})$ such that $\hat{f}=f$ on
$A\times G$ and
\begin{equation*}
\lim\limits_{z\to z_0,\ w\to \eta_0,\
w\in\mathcal{A}_{\alpha}(\eta_0)}\hat{f}(z,w)=f(z_0,\eta_0),\qquad 0<\alpha<\frac{\pi}{2},
\end{equation*}
for every $z_0\in D$ and $\eta_0\in  B\cap B^{\ast}.$ Moreover,   $\vert
\hat{f}\vert_{\widehat{W}^{\text{o}}} \leq\vert f\vert_{W} .$
\end{thm}
\begin{proof} First one proves the existence and  uniqueness of $\hat{f}.$  
Fix an $f :\ W\longrightarrow\C$ which satisfies (i)--(iii) above.

\smallskip

\noindent{\bf Step I:} {\it Reduction to the case where $D\Subset X$ is an open  hyperconvex set\footnote{
An open set $D\subset X$  is said to be {\it hyperconvex} if it admits
an exhaustion function which is bounded subharmonic.} and $A$
is a locally regular  compact subset of $D$
and  $\vert f\vert_W<\infty.$}

\smallskip

Since $X$ is countable at infinity, we may
find an exhaustion sequence  $(D_k)_{k=1}^{\infty}$ of relatively compact,
hyperconvex open subsets $D_k$ of $D$ with $A_k\subset D_k\nearrow D$
(for example, we can choose open subsets  $D_k$ of $D$  with smooth boundary which contains $A_k$).
Similarly, since  $Y$  is countable at infinity,  we may find a sequence $(\Omega_k)_{k=1}^{\infty}$ of open subsets
of $\Omega$ and a sequence  $(B_k)_{k=1}^{\infty}$  of Jordan  measurable subsets of $B$  which satisfy the hypothesis
of Proposition  \ref{prop5.16}. Let  $G_k:=\left\lbrace
   w\in \Omega_k:\ \omega(w,B_k,\Omega_k)<1-\delta\right\rbrace.$
Using a compactness argument,  we see that  $\vert f\vert_{\X(A_k,B_k;D_k,G_k)}<\infty.$

By reduction assumption, for each $k$ there exists an
$\hat{f}_k\in
\mathcal{O}\left(\widehat{\X}^{\text{o}}(A_k,B_k; D_k,G_k)\right)$
such that $\hat{f}_k$ admit the angular limit $f|_{\X(A_k,B_k\cap B^{\ast}_k; D_k,G_k)}$ on $\X(A_k,B_k\cap B^{\ast}_k; D_k,G_k).$

 We  claim that $\hat{f}_{k+1}=\hat{f}_k$ on   $\widehat{\X}^{\text{o}}(A_k,B_k; D_k,G_k).$
Indeed, fix an arbitrary  $k_0\geq 1$ and an arbitrary point $(z_0,w_0)
\in\widehat{\X}^{\text{o}}(A_{k_0},B_{k_0}; D_{k_0},G_{k_0}).$ Let $k\in\N$ such that
$k\geq k_0.$ Let $\mathcal{D}$ be the connected component containing $z_0$
of the following open set
\begin{equation*}
     \left\lbrace z\in  D:\ \omega(z,A_{k_0},D_{k_0}) <1 - \omega_{\delta}
     (w_0,B_k,\Omega_k)\right\rbrace.
     \end{equation*}
Observe that both functions $\hat{f}_{k_0}( \cdot,w_0)|_{\mathcal{D}}$
and $\hat{f}_{k}(\cdot,w_0)|_{\mathcal{D}}$ are holomorphic and
\begin{equation*}
 \hat{f}_k(z,w_0)= f_k(z,w_0)=\hat{f}_{k_0}(z,w_0),\qquad z\in A_k\cap
 \mathcal{D}.
     \end{equation*}
Since  $A_k\cap
 \mathcal{D}$ is non-polar,
  we deduce that  $\hat{f}_{k_0}(\cdot,w_0)|_{\mathcal{D}}=
\hat{f}_{k}(\cdot,w_0)|_{\mathcal{D}}.$ Hence, $\hat{f}_{k_0}(z_0,w_0)=
\hat{f}_{k}(z_0,w_0),$ which proves the above assertion.

On the other hand, by Proposition  \ref{prop5.16} one gets $\widehat{\X}^{\text{o}}(A_k,B_k; D_k,G_k) \nearrow \widehat{W}^{\text{o}} $
as $k\nearrow \infty.$
Therefore, we may glue $\hat{f}_k$ together to obtain
 a function
$\hat{f}\in  \mathcal{O}(\widehat{X}^{\text{o}})$
such that $\hat{f} $ admits the angular limit $f$ on  $ W $ and
$\hat{f}=f$ on $A\times G.$
The uniqueness of such an extension $\hat{f}$  can be proved using the
argument given in the previous paragraph.

This completes Step I.

\smallskip

\noindent{\bf Step II:} {\it The case where
   $D\Subset X$ is an  open hyperconvex set,  $A$ is a locally
 regular  compact subset of $D,$
and  $\vert f\vert_W<\infty.$}

Suppose without loss of generality that $\vert f\vert_W <1.$
We will apply Th\'eor\`eme 3.3 in the work of Zeriahi \cite{zr}  to the pair of condenser
$(A,D).$
%
In the sequel, we will use the notation from this work.

Let $\mu_0:=\mu_{A,D}$  and $\mu_1$  a $B$-admissible Lebesgue  measure of $D.$  Let
    $H_1:= L^2_{h}(D,\mu_1),$ $H_0:=$ the closure of
$H_1|_{A}$ in $L^2(A,\mu_0),  $  let $(b_j)_{j=1}^{\infty}\subset H_1$ be  a   system of
 doubly orthogonal bases in $H_1$ and $H_0.$  Recall that $\Vert b_j\Vert_{H_0}=1.$   Putting
  $\gamma_j:=\Vert b_j\Vert_{H_1},\ j\in\N,$
  we have that
\begin{equation}\label{eq8.2.0}
\sum\limits_{j=1}^{\infty}\gamma_j^{-\epsilon}<\infty,\qquad \epsilon>0.
\end{equation}
For any $w\in B,$ we have $f(\cdot,w)\in H_1$ and   $f(\cdot,w)|_{A}\in
H_0.$ Hence
\begin{equation}\label{eq8.2.1}
f(\cdot,w)=\sum\limits_{j=1}^{\infty} c_j(w)b_j,
\end{equation}
where
\begin{equation}\label{eq8.2.2}
c_j(w)=\frac{1}{\gamma_j^2}\int\limits_{D} f(z,w)\overline{b_j(z)}
d\mu_1 z=\int\limits_{A} f(z,w)\overline{b_j(z)}
d\mu_0(z),\qquad j\in\N.
\end{equation}
Taking  the hypotheses (i)--(iii)
into account and applying    Lebesgue's Dominated Convergence
Theorem, we see that the formula
\begin{equation}\label{eq8.2.3}
\widehat{c_j}(w):=\int\limits_{A} f(z,w)\overline{b_j(z)}
d\mu_0(z),\qquad w\in G\cup B,\  j\in\N;
\end{equation}
defines a bounded  function which is holomorphic in $G.$
Moreover, by (iii) and (\ref{eq8.2.2})--(\ref{eq8.2.3})   it follows that
\begin{equation}\label{eq8.2.4}
\lim_{ w\to \eta,\ w\in \mathcal{A}_{\alpha}(\eta) }\widehat{c_j}(w)  =\widehat{c_j}(\eta)=c_j(\eta),\qquad
 \eta\in B,\ 0<\alpha<\frac{\pi}{2} .
\end{equation}

  Using (\ref{eq8.2.2})--(\ref{eq8.2.4}),  we  obtain the following estimates
\begin{eqnarray*}
\frac{\log{\vert \widehat{c_j}(w)\vert} }{\log{\gamma_j}}&\leq &
\frac{\log{\sqrt{\mu_0(A)}}
}{\log{\gamma_j}},\qquad w\in G,\ j\in\N,\\
 \limsup_{ \ w\to \eta,\ w\in\mathcal{A}_{\alpha}(\eta) }\frac{\log{\vert \widehat{c_j}(\eta)\vert} }{\log{\gamma_j}}&\leq &
\frac{\log{\sqrt{\mu_1(D)}    } }{\log{\gamma_j}}-1, \qquad \eta\in
B,\ 0<\alpha<\frac{\pi}{2},\
j\in\N.
\end{eqnarray*}

This shows that for any $\epsilon>0,$ there is a sufficiently large $N$ such that
for all $j\geq N,$
\begin{equation}\label{eq8.2.5}
\frac{\log{\vert \widehat{c_j}\vert}
}{\log{\gamma_j}}\leq\omega_{\delta}(\cdot,B,\Omega)+\epsilon-1 \qquad\text{on}\ G.
\end{equation}

Take a compact set $K\Subset D$ and let $1>\alpha=\alpha(K) >\max\limits_{K}
\omega(\cdot,A,D).$ Choose an   $\epsilon=\epsilon(K)>0$ so small   that $\alpha+2\epsilon<1.$
Consider the open set
\begin{equation*}
G_K:=\left\lbrace w\in G:\ \omega_{\delta}(\cdot,B,\Omega)<1-\alpha-2\epsilon
\right\rbrace.
\end{equation*}
By (\ref{eq8.2.5}) there is a constant $C^{'}(K)$ such that
\begin{equation}\label{eq8.2.6}
\vert  \widehat{c_j}\vert_{G_j}\leq C^{'}(K)\gamma_j^{\omega_{\delta}(\cdot,B,\Omega)+\epsilon-1}
\leq  C^{'}(K)\gamma_j^{-\alpha-\epsilon},\qquad j\geq 1.
\end{equation}
Now we wish to show that
\begin{equation}\label{eq8.2.7}
\sum\limits_{j=1}^{\infty}\widehat{c_j}(w)b_j(z)
\end{equation}
converges locally uniformly in $\widehat{W}^{\text{o}}.$
Indeed, by (\ref{eq8.2.0})  and (\ref{eq8.2.6}), we have that
\begin{equation}\label{eq8.2.8}
\sum\limits_{j=1}^{\infty}\vert \widehat{c_j}\vert_{G_K}\vert
b_j\vert_K\leq \sum\limits_{j=1}^{\infty}  C^{'}(K)\gamma_j^{-\alpha-\epsilon}
C(K,\alpha)\gamma_j^{\alpha }
\leq  C^{'}(K)C(K,\alpha)\sum\limits_{j=1}^{\infty}\gamma_j^{-\epsilon
}<\infty,
\end{equation}
which gives the normal convergence on $K\times G_K.$
Since the compact set $K$ and $\epsilon>0$ are arbitrary,
the series in (\ref{eq8.2.7}) converges uniformly on compact subsets of
$\widehat{W}^{\text{o}}.$  Let $\hat{f}$ denote this limit function in
 (\ref{eq8.2.7}).

Fix $z_0\in D$ and $\eta_0\in B\cap B^{\ast}.$ We choose a compact $K_0\Subset D$ so that $K_0$ is a
neighborhood
 of  $z_0.$ Let $\epsilon_0>0.$

 In virtue of (\ref{eq8.2.8}), there is an $N_0$ such that
\begin{equation}\label{eq8.2.9}
\sum\limits_{j=N_0+1}^{\infty}\vert \widehat{c_j}\vert_{G_{K_0}}\vert
b_j\vert_{K_0} <\frac{\epsilon_0}{2}.
\end{equation}
On the other hand, in virtue of (\ref{eq8.2.1})--(\ref{eq8.2.4}), we may
find, for any $0<\alpha<\frac{\pi}{2},$ an open neighborhood $V_{\alpha}$
of $\eta_0$ such that
\begin{equation*}
\left\vert \sum\limits_{j=1}^{N_0}  \widehat{c_j}(w)
b_j (z)- \sum\limits_{j=1}^{N_0}  c_j(\eta_0)
b_j (z) \right\vert<\frac{\epsilon_0}{2},\qquad z\in K_0,\ w\in \mathcal{A}_{\alpha}(\eta_0)\cap V_{\alpha}.
\end{equation*}
This, combined with (\ref{eq8.2.7}) and (\ref{eq8.2.9}), implies that
\begin{equation*}
\limsup\limits_{z\to z_0,\ w\to\eta_0,\ w\in\mathcal{A}_{\alpha}(\eta_0)}
\left\vert \hat{f}(z,w)-f(z_0,\eta_0)\right\vert<\epsilon_0,\qquad
 0<\alpha<\frac{\pi}{2}.
\end{equation*}
Since $\epsilon_0>0$ and $(z_0,\eta_0)\in D\times (B\cap B^{\ast})$  can be arbitrarily
chosen, we conclude that
\begin{equation*}
\lim\limits_{z\to z_0,\ w\to\eta_0,\ w\in\mathcal{A}_{\alpha}(\eta_0)}
 \hat{f}(z,w)=f(z_0,\eta_0),\qquad (z_0,\eta_0)\in  D\times (B\cap B^{\ast}),\
 0<\alpha<\frac{\pi}{2}.
\end{equation*}

To complete Step II, it remains to show that $\hat{f}=f$ on $A\times G.$
To do this, fix an arbitrary $(z_0,w_0)\in A\times G.$
 Let $\mathcal{G}$ be the connected component   of  $G$ containing $w_0$. Recall that
 $ G= \left\lbrace w\in  \Omega:\  \omega
     (w,B,\Omega)<1-\delta  \right\rbrace. $ Then
observe that both functions $\hat{f} (z_0,\cdot)|_{\mathcal{G}}$
and $ f(z_0,\cdot)|_{\mathcal{G}}$ admit the same angular limit
$f$ on $B\cap \mathcal{G}^{\Omega}.$ Consequently, applying Theorem \ref{thm6.4}
yields that $\hat{f} (z_0,\cdot)|_{\mathcal{G}}=
 f (z_0,\cdot)|_{\mathcal{G}}.$ Hence, $\hat{f} (z_0,w_0)=
 f (z_0,w_0),$ which proves the above assertion.

This completes the proof of Step II.

It remains to prove the estimate  $\vert \hat{f}\vert_{\widehat{W}^{\text{o}}}\leq\vert f\vert_W.$
In order to reach a contradiction assume   that there is a point
$z^0\in\widehat{W}^{\text{o}}$ such that $\vert\hat{f}(z^0)\vert>\vert
f\vert_W.$ Put $\alpha:=\hat{f}(z^0)$ and consider the function
\begin{equation}\label{eq8.2.10}
g(z):=\frac{1}{f(z)-\alpha},\qquad  z\in W.
\end{equation}
Using the above assumption, it can be checked   that  $g $ satisfies  hypotheses (i)--(iii)
of Theorem \ref{mixedcrossthm}. Hence  applying the first assertion of the theorem, there is exactly
one function  $\hat{g}\in\mathcal{O}(\widehat{W}^{\text{o}} )$ with $\hat{g}=g$ on
$A\times G.$ Therefore, by (\ref{eq8.2.10}) we have on $A\times G:$ $g(f-\alpha)\equiv 1.$ Thus
$\hat{g}(\hat{f}-\alpha)\equiv 1$ on $\widehat{W}^{\text{o}}.$ In particular,
\begin{equation*}
0=\hat{g}(z^0 )(\hat{f}(z^0)-\alpha)= 1;
\end{equation*}
  a contradiction.
Hence the inequality $\vert \hat{f}\vert_{\widehat{W}^{\text{o}}}\leq\vert f\vert_W$ is
proved.
\end{proof}

Finally, we conclude this section with  two uniqueness  results.
\begin{prop}  \label{prop8.5}
  Let $X,\ Y$  be two complex manifolds of dimension $1,$
   $D\subset X,$ $ G\subset Y$  two open sets  and
  $A\subset \partial D,$  $B\subset\partial G$  two Jordan measurable subsets
  of positive  length.
Let   $\tilde{D}\subset  X$   be an open set, $D\cap \tilde{D}\not=\varnothing,$ and let $ \widetilde{A}
\subset \partial \widetilde{D}$
   be a   Jordan measurable subset of positive measure.
Put
\begin{eqnarray*}
 \widehat{W}^{\text{o}}&:=&\widehat{\X}^{\text{o}}\Big(A, B;D, G\Big),\\
\widehat{\widetilde{W}}^{\text{o}}&:=&\widehat{\X}^{\text{o}}\Big(\widetilde{A}, B ; \widetilde{D},G
 \Big ).
\end{eqnarray*}

Let  $\hat{f}\in\mathcal{O}( \widehat{W}^{\text{o}}) ,$
$\hat{\tilde{f}}\in\mathcal{O}( \widehat{\tilde{W}}^{\text{o}}), $
and $z_0\in D\cap \tilde{D}$ be such that
  both  $\hat{f} $  and  $\hat{\tilde{f}} $ admit the same angular
limit  at $\left(z_0,b \right)$ for a.e.
   $ b \in B  . $
Then $\hat{f} (z,w)=\hat{\tilde{f}}(z,w) $
for every $ (z,w)\in  \widehat{W}^{\text{o}}\cap
\widehat{\widetilde{W}}^{\text{o}}.$
\end{prop}
 \begin{proof}
 Fix an arbitrary $w_0\in G $ such that
 $(z_0, w_0)\in\widehat{ W
 }^{\text{o}}\cap\widehat{\widetilde{W}}^{\text{o}}.$
Choose $0<\epsilon<1$ so that
\begin{equation*}
(z_0,w_0)\in D_{1-\epsilon}\times G_{\epsilon}\cap \widetilde{D}_{1-\epsilon}\times G_{\epsilon},
\end{equation*}
where we have  used the notation of level sets introduced in Section 4.
Applying Theorem \ref{thm6.4} to $\hat{f}(z_0,\cdot)|_{ G_{\epsilon}}$ and
 $\hat{\tilde{f}}(z_0,\cdot)|_{ G_{\epsilon}},$
  it follows that
 $
\hat{\tilde{f}}\left(z_0,w_0 \right)=\hat{f}\left(z_0,w_0 \right). $
  Hence,
the proof   is finished.
\end{proof}
Now we are able to  prove the uniqueness stated in Theorem A.
\begin{cor}\label{cor8.6}
We keep the hypotheses and the notation of Theorem A.
Then   there is at most one function $\hat{f}\in
\mathcal{O}(\widehat{W}^{\text{o}})$ which satisfies Property 1) of
Theorem A.
\end{cor}
\begin{proof} It follows immediately from Proposition \ref{prop8.5}.
  \end{proof}
\section{Proof of Theorem A}

 Recall
 that by  Corollary \ref{cor8.6}, the function
$\hat{f}$ satisfying Part 1) is uniquely determined (if it exists).
We only gives here the proof of Part 1).
Using this part, we conclude the proof of Parts 2)--5) of Theorem A in exactly the
same way as we did in Section 6 starting from Step 2 of that section.
The proof is divided into two steps.

\smallskip

\noindent{\bf Step 1:} {\it
 Proof of Theorem A for the case where  $G$
 is a  Jordan   domain.}

\smallskip

 \noindent{\it Proof of Step 1.}  In virtue of Proposition \ref{prop8.2},
 let $\{a_{j}\}_{j\in J}$ be a finite or countable
subset of $A$ with the following properties:
\begin{itemize}
\item[$\bullet$]
 For any $j\in J,$  there is an
open neighborhood $U_{j}$ of  $a_{j}$ such that $D\cap U_{j}$ is   either a
 Jordan domain or the disjoint union of two  Jordan
domains (according to the type of $a_{j}$);
\item[$\bullet$] $A\subset \bigcup\limits_{j\in J} U_{j}.$
\end{itemize}
For any $0<\delta<\frac{1}{2},$   define
\begin{equation*}
\begin{split}
U_{j,\delta}&:=\left\lbrace z\in  D\cap U_{j}:\ \omega(z, A^{\ast}\cap U_{j}, D\cap U_{j})<\delta
 \right\rbrace,\qquad
j\in J,\\
A_{\delta}&:=\bigcup\limits_{j\in J} U_{j,\delta},\\
G_{\delta}&:=\left\lbrace w\in G:\ \omega(w, B,  G)<1-\delta  \right\rbrace.
\end{split}
\end{equation*}

Moreover, for every $j\in J$ let
\begin{equation}\label{eq9.1}
\begin{split}
W_j&:=\X\left(\partial  (D\cap U_{j}) \cap A,B; D\cap
U_{j},G\right),\\
\widehat{W_j}^{\text{o}} &:=\widehat{\X}^{\text{o}}\left(\partial  (D\cap U_{j}) \cap A,B; D\cap
U_{j},G\right),\\
\tilde{f}_j&:=f|_{W_j}.
\end{split}
\end{equation}
Using the hypotheses on $f,$ we conclude that $\tilde{f}_j,$  $j\in J,$
satisfies (i)--(iii) of    Theorem A.  Moreover, since  $G$ is a
Jordan domain and $D\cap U_{j},$  $j\in J,$ is either a Jordan
domain
or the disjoint union of two  Jordan domains, we are able to
apply the result of Section 6  to $\tilde{f}_j.$
Consequently, we obtain, for $j\in J,$  a unique
function $\hat{f}_{j} \in   \mathcal{O}\left(\widehat{W_j}^{\text{o}}\right), $
a subset $A_{j}$ of $\partial (D\cap U) \cap A,$ a subset $B_{j}$ of $B$ such that
\begin{equation}\label{eq9.2}
\begin{split}
&\  A_j\subset  A_j^{\ast},\\
&\ \Big (\partial (D\cap U) \cap A\Big) \setminus A_{j}\ \text{and}\ B\setminus B_{j}\  \text{is of zero length,}\\
\hat{f}_{j}\  &\ \text{admits the angular limit}\ f \  \text{on}\
 \left( \left(\partial  (D\cap U_{j}) \cap A_{j}\right)\times G\right)\cup \left( D\times  B_{j}\right).
\end{split}
\end{equation}
Put
\begin{equation}\label{eq9.3}\begin{split}
\tilde{A}&:=\bigcap\limits_{j\in J} A_{j}\qquad\text{and}\qquad \tilde{B}:=\bigcap\limits_{j\in J}
B_{j},\\
W_{\delta}&:=\X\left( A_{\delta},\tilde{B};D,G_{\delta}\right),\\
\widehat{W_{\delta}}^{\text{o}}&:=\widehat{\X}^{\text{o}}\left( A_{\delta},\tilde{B};D,G_{\delta}\right).
\end{split}
\end{equation}
In virtue of Proposition \ref{prop8.5}, we are able to collect the family $\left(\hat{f}_{j}|_{U_{j,\delta}
\times G_{\delta}} \right)_{j\in J}$ in order to obtain a
  function $\tilde{\tilde{f}}_{\delta}\in \mathcal{O}(A_{\delta}\times
G_{\delta}).$

 Next, consider the function $ \tilde{f}_{\delta}:\ W_{\delta}\longrightarrow
 \C$ given by
\begin{equation}\label{eq9.4}
  \tilde{f}_{\delta}:=
\begin{cases}
 \tilde{\tilde{f}}_{\delta},
  & \qquad\text{on}\  A_{\delta}\times G_{\delta}\\
  f, &   \qquad\text{on}\ D\times (\tilde{B}\cap \tilde{B}^{\ast})
\end{cases} .
\end{equation}
In virtue of (\ref{eq9.1})--(\ref{eq9.4}), we deduce that
\begin{equation}\label{eq9.5}
A\setminus \tilde{A}\ \text{and}\ B\setminus \tilde{B}\  \text{is of zero length,}
\end{equation}
and
\begin{equation}\label{eq9.6}
\begin{split}
\lim\limits_{z\to z_0,\
 w\to b_0,\ w\in\mathcal{A}_{\alpha}(b_0)}\tilde{f}_{\delta}(z,w)&=f(z_0,b_0),
 \qquad 0<\alpha<\frac{\pi}{2},\ z_0\in D,\ b_0\in \tilde{B}\cap\tilde{B}^{\ast},\\
\lim\limits_{z\to a_0,\ z\in\mathcal{A}_{\alpha}(a_0),\
 w\to w_0}\tilde{f}_{\delta}(z,w)&=f(a_0,w_0),
 \qquad 0<\alpha<\frac{\pi}{2},\ a_0\in \tilde{A},\ w_0\in
 G_{\delta}.
\end{split}
\end{equation}
In virtue of (\ref{eq9.4})--(\ref{eq9.6}), $\tilde{f}_{\delta}$ satisfies
the hypotheses (i)--(iii) of Theorem \ref{mixedcrossthm}. Applying this
theorem to $\tilde{f}_{\delta},$ we obtain, for every $0<\delta<\frac{1}{2},$  a
function $\hat{f}_{\delta}\in \mathcal{O}\left(\widehat{W_{\delta}}^{\text{o}}\right).$
In virtue of (\ref{eq9.6}), we see that
\begin{equation}\label{eq9.7}
\begin{split}
\hat{f}_{\delta}&=\tilde{f}_{\delta}\qquad \text{on}\ A_{\delta}\times
G_{\delta},\\
\lim\limits_{z\to z_0,\
 w\to b_0,\ w\in\mathcal{A}_{\alpha}(b_0)}\hat{f}_{\delta}(z,w)&=f(z_0,b_0),
 \qquad 0<\alpha<\frac{\pi}{2},\ z_0\in D,\ b_0\in \tilde{B}\cap \tilde{B}^{\ast},\\
\lim\limits_{z\to a_0,\ z\in\mathcal{A}_{\alpha}(a_0),\
 w\to b_0}\hat{f}_{\delta}(z,w)&=f(a_0,w_0),
 \qquad 0<\alpha<\frac{\pi}{2},\ a_0\in \tilde{A},\ w_0\in
 G_{\delta}.
\end{split}
\end{equation}

We are now in a position to define the desired extension  function $\hat{f}.$
Indeed,
one  glues
$\left(\hat{f}_{\delta}\right)_{0<\delta<\frac{1}{2}}$ together to obtain
$\hat{f}$ in the following way
\begin{equation}\label{eq9.8}
\hat{f}:=\lim\limits_{\delta\to 0} \hat{f}_{\delta}\qquad \text{on}\
\widehat{W}^{\text{o}}=\widehat{\X}^{\text{o}}\left(A,B;D,G\right) .
\end{equation}
Now one has to check that the limit (\ref{eq9.8}) exists and possesses all the required
properties. This will be an immediate consequence of the following
\begin{lem}\label{lem9.1}
 For  any point $(z,w)\in \widehat{W}^{\text{o}} $ put
 \begin{equation}\label{eq9.9}
\delta_{(z,w)} :=\frac{1-\omega(z,A^{\ast},D)-\omega(w,B^{\ast},G)}{2}.
\end{equation}
Then $\hat{f}(z,w)=\hat{f}_{\delta}(z,w)$ for all $0<\delta\leq\delta_{(z,w)}.$
\end{lem}

\smallskip

\noindent{\it Proof of Lemma \ref{lem9.1}.}
Fix an arbitrary point  $(z_0,w_0) \in \widehat{\X}^{\text{o}}\left(A,B;D,G\right)$
and let
$\delta_0:=\delta_{(z_0,w_0)}.$ Let $0<\delta\leq\delta_0.$
Then,  $\omega(w_0,B^{\ast},G)<1-\delta_0$ and
\begin{eqnarray*}
 \omega(z_0,A_{\delta},D)+\omega_{\delta_0}(w_0, B, G)&\leq&
\omega(z_0,A^{\ast},D)+\frac{\omega(w_0, B^{\ast}, G )}{1-\delta_0}\\
&\leq& \frac{\omega(z_0,A^{\ast},D)+\omega(w_0, B^{\ast}, G  )}{1-\delta_0}
<1,
\end{eqnarray*}
where the latter estimate follows from formula (\ref{eq9.9}).
Consequently,
\begin{equation}\label{eq9.10}
(z_0,w_0)\in \widehat{\X}^{\text{o}}\left(A_{\delta} ,B;D,G_{\delta_0}\right).
\end{equation}
 On the other hand, using Part 1) of Proposition \ref{prop8.2}, it is
 clear that
\begin{equation}\label{eq9.11}
 \widehat{\X}^{\text{o}}\left(A_{\delta} ,B;D,G_{\delta_0}\right)
 \subset \widehat{\X}^{\text{o}}\left(A_{\delta} , B;D,G_{\delta}\right)
 \cap \widehat{\X}^{\text{o}}\left(A_{\delta_0} ,B;D,G_{\delta_0}\right) .
\end{equation}
 Moreover,    in virtue of  (\ref{eq9.4}) and (\ref{eq9.7}), we
 have
\begin{equation}\label{eq9.12}
\hat{f}_{\delta} =\tilde{\tilde{f}}_{\delta} = \hat{f}_{\delta_0}  \qquad\text{on}\  A_{\delta}\times
 G_{\delta_0}.
\end{equation}
Next, let  $\mathcal{D}$ be the connected component containing $z_0$ of the
following open set
\begin{equation*}
\left\lbrace z\in D:\
\omega(z,A_{\delta},D)<1-\omega_{\delta_0}(w_0,B,G)\right\rbrace
\end{equation*}
Observe that, in virtue of (\ref{eq9.10})--(\ref{eq9.11}), both functions
$\hat{f}_{\delta}|_{\mathcal{D}}$ and  $\hat{f}_{\delta_0}|_{\mathcal{D}}$
are holomorphic and  $\mathcal{D}\cap A_{\delta}$ is  a nonempty open
set. Therefore,  we deduce from (\ref{eq9.12}) that
 $\hat{f}_{\delta} =\hat{f}_{\delta_0} $ on $\mathcal{D}.$ Hence,
 $\hat{f}_{\delta}(z_0,w_0)=\hat{f}_{\delta_0}
 (z_0,w_0),$
which completes the proof of the lemma. \hfill
$\square$

\smallskip

We complete the proof (of Part 1)) as follows.
An immediate consequence of Lemma \ref{lem9.1} is that
$\hat{f}\in\mathcal{O}\left(\widehat{W}^{\text{o}}  \right).$
 Next, we apply  Lemma \ref{lem9.1} and make use of   (\ref{eq9.4})--(\ref{eq9.9})
  and of  the fact that $\widehat{W_{\delta}}^{\text{o}}\to \widehat{W}^{\text{o}}    $ as
 $\delta\searrow 0.$ Consequently, we conclude   that  $\hat{f} $
 satisfies  the conclusion of Part 1).  Hence, the proof of Step 1 is finished.
 \hfill $\square$

\smallskip

\noindent{\bf Step 2:} {\it Proof of Theorem A for the general case.}

\smallskip

 \noindent{\it Proof of Step 2.}
We proceed using Step 1 in exactly the same way as we proved Step 1
using the result of Section 6.
  Hence, Step 2 is finished.
\hfill  $\square$

  This completes the proof of Theorem A.
\hfill $\square$

 \smallskip

 We conclude this section with the following remark.
  Using the above proof,   one can also derive  Gonchar's Theorem (Theorem 1)  from  Druzkowski's Theorem
   (Theorem 3).
Indeed,  in Step 1 above,
  let $\{a_j\}_{j\in J}$ be  finite or countable subset of $ A$ with the following properties:
  \begin{itemize}
\item[$\bullet$]
 For any $j\in J,$  there is an
open neighborhood $U_{j}$ of  $a_{j}$ such that $D\cap U_{j}$ is    a
 Jordan domain    and  $A\cap U_{j}$  is one open arc;
\item[$\bullet$] $A\subset \bigcup\limits_{j\in J} U_{j}.$
\end{itemize}
    Then we repeat Step 1 ( $B$ is only one open arc) and Step 2 (the general case) above using Druzkowski's Theorem.
  Gonchar's Theorem follows.  \hfill  $\square$

\section{Proof of Theorem B}
  We will only give the proof of Theorem B for the case when
  $D$  and  $G$ are the unit disc $E.$ Since the general case
  can be proved using the scheme of Section 6 and  8,
    it is left to the interested reader.
  The proof is divided into the following  two steps.

  \smallskip

\noindent{\bf Step 1:} {\it Proof of Theorem B for the case when
 the slice functions  $f(a,\cdot)|_{G}$ and $f(\cdot,b)|_{D}$ are bounded
for every $a\in A$ and $b\in B.$}

\smallskip

 \noindent{\it Proof of Step 1.} For any $N\in\N$ let
 \begin{equation}\label{eq11.1.1}
A_{N}:=\left\lbrace a\in A:\ \vert f(a,\cdot)\vert_{G}  \leq N \right\rbrace
\ \text{and}\  B_{N}:=\left\lbrace b\in B:\ \vert f(\cdot,b)\vert_{D}  \leq N
\right\rbrace.
 \end{equation}
Using the assumption of Step 1 and (\ref{eq11.1.1}), we obtain
 \begin{equation}\label{eq11.1.2}
 A_{N}\nearrow A\  \text{and} \  B_N\nearrow B\qquad\text{as}\ N\nearrow \infty.
\end{equation}
Now we would like to show that for every $N\in\N,$
 \begin{equation}\label{eq11.1.3}
 \begin{split}
 A_{N}&\quad\text{is a closed subset of}\ A\ \text{and}\ f|_{A_{N}\times
 G}\in \mathcal{C}(A_{N}\times G),\\
  B_{N}&\quad\text{is a closed subset of}\ B\ \text{and}\ f|_{D\times
 B_N}\in \mathcal{C}(D\times B_{N}) .
 \end{split}
\end{equation}
To do this
 fix an arbitrary  $N\in \N$
and     let $(a_{n})_{n=1}^{\infty}$ be
 a sequence in $A_{N}$ such that $\lim\limits_{n\to\infty} a_{n}=a_0\in A_N .$
 Consequently, by hypothesis (i),
 \begin{equation}\label{eq11.1.4}
 \lim\limits_{n\to\infty} f(a_{n},t)=f(a_0,t),\qquad t\in B.
 \end{equation}

 On the other hand, it follows from the assumption $(a_{n})_{n=1}^{\infty} \subset A_{N}$
 and the hypothesis of Step 1 that
 \begin{equation*}
\vert f(a_{n},\cdot)\vert_{G}\leq N\qquad\text{and}\qquad \vert
f(a_0,\cdot)\vert_{G}<\infty.
 \end{equation*}
Combining this and (\ref{eq11.1.4}), we are able to apply Khinchin--Ostrowski Theorem  (see \cite[Theorem 4, p. 397]
{go}) to the sequence
 $\left( f(a_{n},\cdot)|_{G}\right)_{n=1}^{\infty}\subset \mathcal{O}(G).$
 Consequently, this sequence  converges uniformly on compact subsets of $G$ to $f(a_0,\cdot).$
    This completes the proof of (\ref{eq11.1.3}).

 On the other hand, by hypothesis (ii), the
holomorphic function $f(a,\cdot)$ admits the angular limit $f(a,b)$
at $b\in B.$  Hence,
it follows that
$f|_{A_{N}\times B_{N}}$ is measurable. Moreover, by (\ref{eq11.1.1}),
 $|f|_{\X(A_{N}, B_{N};D,G)}\leq N$ for every $N\in \N.$ In addition,
in virtue of (\ref{eq11.1.2}),  there exists a sufficiently large integer  $N_0$ such that
 $\mes(A_{N})>0$  and $\mes(B_N)>0$ for
 $N\geq N_0.$
Consequently, we are in a position to apply Theorem A to the function
$f$ restricted to the cross
$\X(A_{N},B_{N};D,G)$  for $N\geq N_0.$ Therefore, we obtain a function
$\hat{f}_N\in\mathcal{O}\left(  \widehat{\X}^{\text{o}}\left(A_{N},B_{N};D,G\right)\right)$
and a subset
$\tilde{A}_N$  (resp. $\tilde{B}_N$) of $A_N$  (resp. $B_N$)   for $N\geq N_0,$ such
that
\begin{equation}\label{eq11.1.5}
\begin{split}
&\ \mes(A_{N}\setminus \tilde{A}_{N})=\mes(B_{N}\setminus \tilde{B}_{N})=0,\\
\hat{f}_{N}\  &\ \text{admits the angular limit}\ f \  \text{on}\
\Big(\tilde{A}_{N}\times G\Big)\bigcup\Big(D\times \tilde{B}_{N}  \Big).
\end{split}
\end{equation}
Put
\begin{equation}\label{eq11.1.6}
\tilde{A}:=\bigcup\limits_{N=N_0}^{\infty}\tilde{A}_{N}\qquad\text{and}\qquad
\tilde{B}:=\bigcup\limits_{N=N_0}^{\infty}\tilde{B}_{N}.
\end{equation}
Applying   (\ref{eq11.1.2}), (\ref{eq11.1.5}), and    Corollary
\ref{cor8.6}, we obtain
\begin{equation}\label{eq11.1.8}
\hat{f}_N=\hat{f}_{N+1} \qquad\text{on}\  \widehat{\X}^{\text{o}}\left(\tilde{A}_{N},\tilde{B}_{N};D,G \right)
 ,\ N\geq
N_0.
\end{equation}
Therefore,  we  may glue the $\hat{f}_N$ together to obtain the desired extension function
$\hat{f}$   as
\begin{equation}\label{eq11.1.9}
\hat{f}=\lim\limits_{N\to\infty}\hat{f}_{N} \qquad\text{on}\  \widehat{W}^{\text{o}}:=
\widehat{\X}^{\text{o}}\left(A,B;D,G\right) .
\end{equation}
Moreover, in virtue of  (\ref{eq11.1.5})--(\ref{eq11.1.9}), we get that
\begin{equation}\label{eq11.1.7} \begin{split}
&\ \mes(A\setminus \tilde{A})=\mes(B\setminus \tilde{B})=0,\\
\hat{f}\  &\ \text{admits the angular limit}\ f \  \text{on}\
 \Big(\tilde{A}\times G\Big)\bigcup\Big(D\times \tilde{B}  \Big).
\end{split}
\end{equation}

Next, for every  $N\geq N_0,$ in virtue of (\ref{eq11.1.2})--(\ref{eq11.1.3}) and
(\ref{eq11.1.5}),
one may find   a sequence $(F_{N,n})_{n=1}^{\infty}$
 (resp.  $(H_{N,n})_{n=1}^{\infty}$)
    of compact subsets of
$\partial D$   (resp.  $\partial G$) such that
\begin{equation}\label{eq11.1.10}\begin{split}
& F_{N,n}\subset F_{N,n+1}\subset A,\qquad   H_{N,n}\subset H_{N,n+1}\subset B, \\
 &\mes(F_{N,n})>0,\quad \mes(H_{N,n})>0, \\
 &\mes\left(\tilde{A}_{N}\setminus \bigcup\limits_{n=1}^{\infty}F_{N,n}  \right)=0,\quad
 \mes\left(\tilde{B}_{N}\setminus \bigcup\limits_{n=1}^{\infty}H_{N,n}  \right)=0  .
\end{split}\end{equation}
Moreover, for any  $k\in \N,$ $k \geq 1,$
and for any $m\in\N,$ put
\begin{equation}\label{eq11.1.11}
\begin{split}
A_{Nnmk}&:=\left\lbrace a\in A_{N}:\ \vert f(a,\zeta)-f(a,\eta)\vert \leq \frac{1}{2k^2},
\qquad  \zeta,\eta\in H_{N,n}:\ \vert \zeta-\eta\vert <\frac{1}{m}
\right\rbrace,\\
B_{Nnmk}&:=\left\lbrace b\in B_{N}:\ \vert f(\zeta,b)-f(\eta,b)\vert \leq \frac{1}{2k^2},
\qquad  \zeta,\eta\in F_{N,n}:\ \vert \zeta-\eta\vert <\frac{1}{m}
\right\rbrace.
\end{split}
\end{equation}

Since, by hypothesis (i), $f\in\mathcal{C}_s(A\times B),$ we deduce from (\ref{eq11.1.10})
and   (\ref{eq11.1.11}) that $A_{Nnmk}$  (resp.  $B_{Nnmk}$)
is a closed subset of $A_{N}$  (resp.  $B_N$) and
\begin{equation}\label{eq11.1.12}
A_{Nnmk}\nearrow A_{N}\  \text{and}\  B_{Nnmk}\nearrow B_{N}\    \qquad\text{ as}\  m\nearrow \infty,\ k\geq
1.
\end{equation}
 Consequently, there is an $m_0:=m_0(N,n,k)$ such that     $\mes(A_{Nnmk}\cap F_{N,n})>0$
 and   $\mes(B_{Nnmk}\cap H_{N,n})>0$ for
 any  $m>m_0.$ Now we are in a position to apply Theorem A to the function $f$ restricted on the cross
$\X\left( A_{Nnmk}\cap F_{N,n},B_{Nnmk}\cap H_{N,n};D,G\right).$
Using (\ref{eq11.1.8})--(\ref{eq11.1.7}) and  Corollary \ref{cor8.6}, we
obtain exactly the   function $\hat{f}$ restricted to
$\widehat{\X}^{\text{o}}\left( A_{Nnmk}\cap F_{N,n},B_{Nnmk}\cap H_{N,n};D,G\right).$
Let\footnote{Recall from Subsection 2.2 that for a boundary subset $T,$  $T^{\ast}$ denotes as usual the set of locally regular
 points relative to $T.$}

\begin{equation}\label{eq11.1.13}
\begin{split}
\tilde{A}_{Nnmk}&:= \left(A_{Nnmk}\cap F_{M,n}\right)\cap\left( A_{Nnmk}\cap
F_{N,n}\right)^{\ast},
\\
\tilde{B}_{Nnmk}&:= \left(B_{Nnmk}\cap H_{M,n}\right)\cap\left( B_{Nnmk}\cap
H_{N,n}\right)^{\ast},
\end{split}
\end{equation}
Taking (\ref{eq11.1.11})--(\ref{eq11.1.13}) into account and arguing as in Step 5 of Section 6, we
may show that
\begin{equation}\label{eq11.1.14}
\begin{split}
\mes\left(\tilde{A}_{Nnmk}\setminus F_{N,n}
\right)=0,\quad   \mes\left(\tilde{B}_{Nnmk}\setminus H_{N,n}
\right)&=0, \\
\limsup\limits_{ (z,w)\to  (a,b)\
z\in\mathcal{A}_{\alpha}(a),  \  w\in\mathcal{A}_{\alpha}(b),
 \ (z,w)\in \widehat{X}^{\text{o}}  }\vert\hat{f}(z,w)
-f(a,b)\vert &<\frac{1}{k},\qquad 0<\alpha<\frac{\pi}{2},
\end{split}
\end{equation}
for every $(a,b)\in \tilde{A}_{Nnmk}  \times \tilde{B}_{Nnmk}.$
Now it suffices to  put
\begin{equation*}
\tilde{A}:=\bigcap\limits_{k=1}^{\infty}  \bigcup\limits_{N=N_0}^{\infty}
  \bigcup\limits_{n=1}^{\infty} \bigcup\limits_{m=m_0(N,n,k)}^{\infty}  \tilde{A}_{Nnmk}\  \text{and}\
\tilde{B}:=\bigcap\limits_{k=1}^{\infty}  \bigcup\limits_{N=N_0}^{\infty}
  \bigcup\limits_{n=1}^{\infty} \bigcup\limits_{m=m_0(N,n,k)}^{\infty}  \tilde{B}_{Nnmk}.
\end{equation*}
Combining this and  (\ref{eq11.1.14}),  (\ref{eq11.1.12}), (\ref{eq11.1.7}) and (\ref{eq11.1.2}),
 we may check that all the conclusions  of Theorem B are
satisfied.  Hence the proof is complete in this first step.
\hfill $\square$

\smallskip

\noindent{\bf Step 2:} {\it  The general case.}

\smallskip

 \noindent{\it Proof of Step 2.}
  We begin with the following
\begin{defi} \label{defi11.1}
For  a closed subset $F$ of $\partial E$ and an   $n\in \N$  with $n>1,$ define the following
 open set
\begin{equation*}
\Delta=\Delta(F,n):=\bigcup\limits_{\zeta\in F}
\left\lbrace z\in \mathcal{A}_{\frac{\pi}{4}}(\zeta):\ \vert z\vert\geq 1-\frac{1}{n}  \right\rbrace
\cup \B\left(0,1-\frac{1}{n}\right).
\end{equation*}
\end{defi}
The reader should compare this definition with Definition \ref{def6.1}.
Below we give   a list of properties of such    open sets.
\begin{prop} \label{prop11.2}  Let $F$ be a closed subset  of $\partial E.$ \\
 1) Let $ \Delta(F,n)$ be   as in
Definition \ref{defi11.1}, then $\Delta(F,n)$ is a rectifiable Jordan domain
and $F\subset \partial\Delta(F,n).$\\
2)  $\Delta(F,n)\nearrow E$ as $n\nearrow\infty.$\\
3) Consider a locally bounded function $f: \ E\cup F\longrightarrow \C.$
Then $\vert f\vert_{\Delta(F,n)}<\infty$ for every $n\in\N$ with $n>1.$\\
4) There holds the following equality
\begin{equation*}
\omega(z,F,E)=\lim\limits_{n\to\infty}\omega\left(z,F,\Delta(F,n)\right),\qquad
z\in E.
\end{equation*}
\end{prop}

\smallskip

\noindent{\it Proof of Proposition \ref{prop11.2}.}
Part 1) may be done   as in the proof of Proposition \ref{prop6.2}.

Part 2) is an immediate consequence of Definition \ref{defi11.1}.

Part 3) follows immediately from the compactness of $F.$

The proof of Proposition \ref{prop5.16} still works in the context of Part
4) making the obviously necessary changes. This completes Part 4).
\hfill $\square$

\smallskip

Now we are in a position to complete Step 2. Indeed, first suppose that
both $A$ and $B$ are closed. Then  consider the   sequence of
rectifiable Jordan domain  $(D_{n})_{n=2}^{\infty}$ and   $(G_{n})_{n=2}^{\infty}$  given by
\begin{equation*}
 D_{n}:=\Delta(A,n)\  \text{and}\  G_n:=\Delta(B,n),\qquad n\in \N,\ n>1.
\end{equation*}
For $n\in\N,$ $n>1,$ let $f_n:=f|_{\X(A,B; D_{n},G_{n})}.$ In virtue of Proposition
\ref{prop11.2}, we are able to apply the result of Step 1 to $f_n.$
Consequently, we obtain a function $\hat{f}_n\in \widehat{\X}^{\text{o}}
(A,B; D_{n},G_{n}).$
Therefore,  we  may glue $\hat{f}_n$ together in order  to  obtain the desired extension function
$\hat{f}$   as
\begin{equation*}
\hat{f}=\lim\limits_{n\to\infty}\hat{f}_{n} \qquad\text{on}\  \widehat{W}^{\text{o}}=
\widehat{\X}^{\text{o}}\left(A,B;D,G\right) .
\end{equation*}
Using Proposition \ref{prop11.2}, we can show that $\hat{f}$
possesses all the   assertions of Theorem B.

The case when $A$ and $B$ are only measurable is similar.
It suffices to  find a sequence $(A_{m})_{m=1}^{\infty}$ of subsets of
$A$ such that $A_{m}$ is compact and $\mes\left( A\setminus
\bigcup\limits_{m=1}^{\infty} A_{m}\right)=0,$  and a similar sequence   $(B_{m})_{m=1}^{\infty}$  for $B.$
 Then we  may apply the
previous discussion to $f|_{\X(A_{m}, B_{m};D,G)}$ in order to
 obtain a function $\hat{f}_m\in \widehat{\X}^{\text{o}}
(A_{m},B_{m}; D,G),$  and define the
desired extension function
$\hat{f}$ by $
\hat{f}:=\lim\limits_{m\to\infty}\hat{f}_{m}$ on $  \widehat{W}^{\text{o}}.$
This completes the proof in
this last step.
\hfill  $\square$
\section{Examples and Concluding remarks}
The following  examples of Dru\.{z}kowski \cite{dr} show the optimality of Theorem A
and B.

Consider   $D=G=E,$  $A=B=
\left\lbrace t\in\partial E: \Re t>0 \right\rbrace,$
$W:=\X(A,B;D,G),$ and $T:=(D\cup A)\times (G\cup B).$

\noindent{\bf Example 1.} Define a function $h:\ T\longrightarrow \C$ as
follows
\begin{equation*}
 h(z,w):=
\begin{cases}
\exp\Big(-\left\lbrack  \Log{(1-z)}+\Log{(1-w)}\right\rbrack \Log{\frac{2+zw}{3}}\Big),
  & z\not=1,\ w\not =1 \\
 0, &  z=1\ \text{or}\ w =1
\end{cases}.
\end{equation*}
where $\Log$ is the principal branch of logarithm.

Put $f:=h|_W.$ As in \cite{dr} observe  that $ f$ is measurable, $f\in \mathcal{C}_s(W)\cap
\mathcal{O}_s(W^{\text{o}}),$ $\vert f\vert_W<\infty,$ but $f|_{A\times B}$ is not continuous
at $(1,1).$  Since $h|_{\widehat{W}^{\text{o}}}\in \mathcal{O}(\widehat{W}^{\text{o}}),$ using the uniqueness
established in Theorem A, we conclude that the solution $\hat{f}$ provided by Theorem A and B
 satisfies $\hat{f}=h|_{\widehat{W}^{\text{o}}}.$
 In addition, we see that, for $0<\alpha<\frac{\pi}{2},$ the
angular limit of $ \hat{f} $ at $(1,1)$ does not
exist. Thus the condition in assertion 3) of Theorem A is necessary.
Moreover, the sets $\tilde{A},$ $\tilde{B}$ given by Theorem B do depend on
$f.$

\noindent{\bf Example 2.} Define a function $h:\  T\longrightarrow \C$ as
follows
\begin{equation*}
 h(z,w):=
\begin{cases}
\exp\left(- ( z-\lambda)\Log^2{\frac{3+w}{1-w}}\right),
  &   w\not =1 \\
 0, &   w =1
\end{cases}.
\end{equation*}
where $(z,w)\in T,$ $0<\lambda\leq \frac{\sqrt{2}}{2}.$

Define $f:=h|_W.$ Then $\hat{f}=h|_{\widehat{W}^{\text{o}}}.$  As in \cite{dr} observe  that $f|_{A\times B}$ is  continuous,
 $f\in \mathcal{C}_s(W)\cap
\mathcal{O}_s(W^{\text{o}}),$ but $ f$ is not locally bounded on $W. $

 In addition, for $\frac{\pi}{3}<\alpha<\frac{\pi}{2},$  consider the functions
 $z_{\alpha,\lambda},w_{\alpha}:\ [0,1]\to \C$
  given by
\begin{eqnarray*}
w_{\alpha}(t)&:=&1+te^{i\left(\pi- \frac{9\alpha}{10}\right)},\\
z_{\alpha,\lambda}(t)&:=& \lambda+\left( \Re\Log^2{\frac{3+w_{\alpha}(t)}{1- w_{\alpha}(t)}}
\right)^{-1}+i\lambda,\qquad t\in [0,1].
\end{eqnarray*}
We may prove that there is an $t_{\alpha,\lambda}>0$ and a neighborhood $U_{\alpha,\lambda}$ of
$\lambda+i\lambda $ in $\C$ such that
\begin{equation*}
\left (z_{\alpha,\lambda}(t),w_{\alpha}(t)\right)\in
\begin{cases}
\Big(\Big(\mathcal{A}_{\alpha}(\lambda+i\lambda)\cap U_{\alpha,\lambda}\Big) \times\mathcal{A}_{\alpha}(1)
\Big)\cap\widehat{W}^{\text{o}} ,
  &  0<t<t_{\alpha,\lambda},\ \lambda =\frac{\sqrt{2}}{2} \\
\Big (U_{\alpha,\lambda}\times\mathcal{A}_{\alpha}(1)\Big )\cap \widehat{W}^{\text{o}} , &   0<t<t_{\alpha,\lambda}
 ,\  0< \lambda <\frac{\sqrt{2}}{2}
\end{cases}.
\end{equation*}
In addition, it can be checked that
\begin{equation*}
\lim_{t\to 0}\left(z_{\alpha,\lambda}(t),
w_{\alpha}(t)\right)=(\lambda+i\lambda,1)\  \qquad\text{and}\  \qquad  \lim_{t\to 0}\left\vert
\hat{f}\left(z_{\alpha,\lambda}(t),w_{\alpha}(t)\right)\right\vert=\infty.
\end{equation*}
This shows that the assumption of the  local boundedness on $f$ is necessary in
Theorem A.

\medskip

Finally, we conclude the article by some remarks and open questions.

\medskip

\noindent 1. It may be proved that $\widehat{W}^{\text{o}}$
provided by Theorem A is the maximal domain of holomorphic
extension of the function $f.$ We postpone the proof of this
result to an ongoing work (see \cite{pn3}).

\smallskip

\smallskip

\noindent 2. Does  Theorem A still hold if we  omit the assumption
(ii) ``$f|_{A\times B}$ is   Jordan-measurable"?

 \smallskip

\noindent 3. Does  Theorem B still hold if we  omit the  assumption  that $f|_{A\times B}\in\mathcal{C}_s(A\times B)$?

\smallskip

\end{document}